\documentclass[12pt]{article}

\usepackage{graphicx}           
\usepackage{epstopdf}
\usepackage{color}              
\usepackage{xcolor,graphicx}
\usepackage{geometry}
\usepackage{float}
\usepackage{times}
\usepackage{indentfirst}        
\usepackage{amsmath,amssymb,bm} 
\usepackage{amsthm}
\usepackage{cases}              
\usepackage{setspace}

 \usepackage[ascii]{inputenc}

\usepackage{titlesec}
\usepackage[all]{xy}            
\usepackage{tikz}
\usetikzlibrary{arrows.meta}
\usetikzlibrary{shapes.gates.ee}
\usetikzlibrary{calc}
\usepackage{multirow}
\usepackage{graphicx} 
\usepackage{epstopdf}
\usepackage{caption}
\usepackage{diagbox}
\usepackage{mathrsfs}
\usepackage{dsfont}
\usepackage[pagebackref]{hyperref}
\hypersetup{ pagebackref, colorlinks=true}
\usetikzlibrary{fit}
\usepackage{tcolorbox}
\usepackage{fancyhdr,fancyvrb}
\usetikzlibrary{shapes.geometric}  

\usepackage{diagbox}
\usepackage{subfigure}
\usepackage{subcaption}
\usepackage{ulem}
\usepackage{booktabs}

\newtheorem{thm}{Theorem}[section]
\newtheorem{defi}[thm]{Definition}

\newtheorem{prop}[thm]{Proposition}
\newtheorem{cor}[thm]{Corollary}

\renewcommand{\arraystretch}{2} 

\makeatletter \@addtoreset{equation}{section} \makeatother
\makeindex 
\def\qed{\hfill \rule{4pt}{7pt}}

\begin{document}

\begin{center}

{\Large \bf  Labeled Plane Trees and  Increasing Plane Trees }

\vskip 4mm
Lora R. Du $^1$, Kathy Q. Ji $^2$ and Dax T.X. Zhang $^3$
\vskip 2mm

$^{1,\, 2}\,$Center for Applied Mathematics, KL-AAGDM\\[2pt]
Tianjin University\\[2pt]
Tianjin 300072, P.R. China \\ \vskip 0.5cm 

$^3\,$College of Mathematical Science \\[2pt]
Institute of Mathematics and Interdisciplinary Sciences \\[2pt]

Tianjin Normal University\\[2pt]
Tianjin 300387, P. R. China

\vskip 2mm
 E-mails: $^1$loradu@tju.edu.cn,   $^2$kathyji@tju.edu.cn and  $^3$zhangtianxing6@tju.edu.cn

\end{center}

\vskip 2mm
\noindent
{\bf Abstract:}  The main aim of this paper is to establish a polynomial analogue of $(n+1)!C_n=2^n(2n-1)!!$ (with $C_n$ as the $n$-th Catalan number)  in the setting of labeled plane trees and increasing plane trees. This analogue is formulated in terms of    improper edges of labeled plane trees and yields explicit formulas for the generating polynomials defined on labeled plane trees refined by improper and proper edges, together with a root-degree  refinement for trees rooted at $0$.  To prove this result, we construct a new involution on labeled plane trees, which implies that  the number of improper edges and the number of proper edges are equidistributed over the set of labeled plane trees. We further apply this involution to establish pairwise symmetry properties of multivariable polynomials defined on labeled plane trees involving several classes of leaves and interior vertices. More precisely, certain specializations of these polynomials are invariant   under the subgroup of $S_6$  generated by the three disjoint transpositions $(12)$, $(34)$, and $(56)$. As special cases, our results  recover the symmetry properties for  plane trees and tip-augmented plane trees due to Dong, Du, Ji and Zhang. 
 Finally, via the Koganov--Janson correspondence,   improper edges of labeled plane trees correspond bijectively to improper arcs of quasi-Stirling permutations, leading to an explicit formula for the generating function defined on quasi-Stirling permutations  refined by improper arcs.

\vskip 2mm



 \vskip 2mm

\section{Introduction}

This paper aims to build a connection between labeled plane trees and increasing plane trees based on the notion of improper edges in labeled plane trees.   Recall that a plane tree is a rooted tree in which the children of each node are linearly ordered. It is well-known that the number of labeled plane trees with $n$ edges is $(n+1)!C_n$, where $C_n=\frac{1}{n+1}{2n \choose n}$ is the $n$-th Catalan number. An increasing plane tree, often referred to as a plane recursive tree, is a plane tree in which the node labels increase along any path from the root. Note that plane recursive trees also appear in literature under the names plane-oriented recursive trees, heap-ordered trees, and sometimes also as scale-free trees. The number of increasing plane trees with $n$ edges is $(2n-1)!!=1\cdot 3 \cdots (2n-1)$.  It is clear that  
\begin{equation}\label{aa}
(n+1)!C_n=2^n(2n-1)!!.
\end{equation}

Our primary  goal is to give a polynomial analogue of \eqref{aa}  based on the notion of  improper edges in labeled plane trees. Improper edges of labeled plane trees  were introduced by  Guo and Zeng \cite{Guo-Zeng-2007}  to offer a combinatorial interpretation of   the generalized Ramanujan polynomials defined by Chapoton \cite{Chapoton-2002}. It is worth noting that Shor \cite{Shor-1995} and Dumont and Ramamonjisoa \cite{Dumont-Ramamonjisoa-1996} independently defined the improper edge of  labeled rooted trees,  which resulted in  a refinement of Cayley's formula.   Zeng \cite{Zeng-1999}   later discovered a link between Shor's polynomials and the Ramanujan polynomials; see  Chen, Fu and Wang \cite{Chen-Fu-Wang-2025} and Chen and Yang \cite{Chen-Yang-2021} for more details.  Subsequently,   Guo and Zeng \cite{Guo-Zeng-2007} extended the notion of  improper edges from labeled rooted trees to labeled plane trees.

  \begin{figure}[H]
    \centering
\begin{tikzpicture}
 [scale=0.8,vertex/.style={shape=circle, draw, inner sep=1.5pt, fill=black},
 subtree/.style={shape=ellipse, draw,minimum width=1.5cm, minimum height=.5cm},
 every fit/.style={ellipse,draw,inner sep=-2pt},
sibling distance=1.5cm,level distance=14mm,
 leaf/.style={label={[name=#1]below:$ $}},auto]
 
\node[vertex,label=0:{$1$}]{}[grow=down, sibling distance=2.5cm,level distance=15mm]
  child {node [vertex, label=-180:{$6$}]{} 
  child {node [vertex, label=-90:{$0$}]{} edge from parent[solid] }edge from parent[solid]}
  child {node [vertex, label=0:{$3$}]{}[sibling distance=10mm,level distance=15mm]
  child{node[vertex,label=-90:{$7$}] {}edge from parent[solid]  } child{node[vertex,label=-90:{$2$}] {}edge from parent[solid]} child{node[vertex,label=-90:{$5$}] {}edge from parent[solid]} child{node[vertex,label=-90:{$4$}]{} edge from parent[solid] }edge from parent[solid]};
\end{tikzpicture}
\caption{A labeled plane tree with $7$ edges.}
\label{An example for labeled plane tree}
    \end{figure}

\begin{figure}[H]
    \centering
 \begin{tikzpicture}
 [scale=0.6,vertex/.style={shape=circle, draw, inner sep=1.5pt, fill=black},
 subtree/.style={shape=ellipse, draw,minimum width=1.5cm, minimum height=.5cm},
 every fit/.style={ellipse,draw,inner sep=-2pt},
sibling distance=2cm,level distance=14mm,
 leaf/.style={label={[name=#1]below:$ $}},auto]
 \node[vertex,name=1,label=0:{$0$}]{}[grow=down]
  child {node [vertex, name=3, label=180:{$1$}]{}  
  child {node [vertex, name=5,label=0:{$2$}]{}
  child {node [vertex, label=-90:{$7$}]{} }
  child {node [vertex, label=0:{$3$}]{}
  child {node [vertex, label=-90:{$5$}]{} }
  child {node [vertex, label=-90:{$4$}]{}  }}} }
  child {node [vertex, label=-90:{$6$}]{}  };
\end{tikzpicture} 
 
\caption{An increasing plane tree with $7$ edges.}
\label{An example for increasing plane tree}
\end{figure}
 
Let $\mathcal{P}_n$ denote the set of labeled plane trees with $n$ edges whose vertices are labeled by $0,1,\ldots,n$, and let $\mathcal{O}_n$ denote the subset of $\mathcal{P}_n$ with  root $0$. Given a vertex $j$ of a labeled plane tree $T$, let $\beta(j)$ be the smallest label in the subtree rooted at $j$. Suppose that the vertex $j$ has children: $j_1,j_2,\ldots,j_k$  (ordered left to right). We call an edge $(j,j_i)$ improper  if  
\begin{equation}
    \beta(j_i)< \min\{ j,\beta(j_{i+1}),\ldots,\beta(j_{k})\},
\end{equation}
otherwise we call it  a proper edge.

For example, the labeled plane tree $T \in \mathcal{P}_{7}$ depicted in Fig.~\ref{An example for labeled plane tree} has three  improper edges: $(1,6), (6,0)$, and $(3,2)$. All other edges in $T$ are proper.
 
For any tree $T$, let ${\rm impr}(T)$ be the number of improper edges of $T$ and let ${\rm prop}(T)$ be the  number of proper edges of $T$. For the labeled plane tree $T \in \mathcal{P}_{7}$ shown in Fig.~\ref{An example for labeled plane tree}, this  gives  
${\rm impr}(T)=3$ and ${\rm prop}(T)=4$.

It is worth noting that every edge in an increasing plane tree is proper. Conversely,  any labeled plane tree without improper edges must be an increasing plane tree.  

For example,  in  the increasing plane tree $T \in \mathcal{P}_{7}$ shown in Fig.~\ref{An example for increasing plane tree}, we have 
${\rm impr}(T)=0$ and ${\rm prop}(T)=7$.

In this paper, we consider the polynomials  that incorporate both improper edges and proper edges of  labeled plane trees: 
\begin{align*}
P_n(x,y)&=\sum_{T\in \mathcal{P}_{n}} x^{{\rm impr}(T)}
 y^{{\rm prop}(T)},\\[5pt]
O_n(x,y,t)&=\sum_{T\in \mathcal{O}_{n}} x^{{\rm impr}(T)} y^{{\rm prop}(T)-{\rm deg}_T(0)} t^{{\rm deg}_T(0)},
\end{align*}
with the convention that $P_0(x,y)=O_0(x,y,t)=1$, where ${\rm deg}_T(0)$ denotes the number of children of root $0$.  Since all edges incident with the root $0$ are proper, these root edges are recorded separately by the variable $t$.
 
 For example, 
$$
\begin{aligned}
    &P_1(x,y)=x+y;\\[5pt]
    &P_2(x,y)=3x^2+6xy+3y^2;\\[5pt]
    &P_3(x,y)=15x^3+45x^2y+45xy^2+15y^3.
\end{aligned}
$$
And 
$$
\begin{aligned}
    &O_1(x,y,t)=t;\\[5pt]
    &O_2(x,y,t)=2t^2+tx+ty ;\\[5pt]
    &O_3(x,y,t)=6t^3+6t^2x+6t^2y+3tx^2+6txy+3ty^2.
\end{aligned}
$$

 By constructing an involution on labeled plane trees, we show the following consequence: 

\begin{thm}\label{main1}  For $n\geq 1$,  
\begin{align}
P_{n}(x,y)&=(2n-1)!! (x+y)^n,\label{Pnxy}\\
O_n(x,y, t) &= \sum_{r=1}^n t^r S_{n,r} (x+y)^{n-r},\label{Fnxyt}
\end{align}
where $S_{n,r}$ counts the number of increasing plane trees with $n$ edges whose root $0$ has degree  $r$. 
\end{thm}

Let 
\begin{equation}
S_n(t)=\sum_{r=1}^n S_{n,r} t^r
\end{equation}
with the convention that $S_0(t)=1$, where  $S_{n,r}$ counts the number of increasing plane trees with $n$ edges so that the degree of $0$ is $r$. 
Evidently,  $S_n(1)=(2n-1)!!$, so we have 
\begin{equation}\label{S(t,q)-gfa}
\sum_{n=0}^\infty  S_n(1) \frac{q^n}{n!}=\frac{1}{\sqrt{1-2q}}.
\end{equation} 
Using the compositional formula \cite[Theorem 5.1.4]{Stanley-2024}, one can readily show that   
\begin{equation}\label{S(t,q)-gf}
 \sum_{n=0}^\infty S_n(t) \frac{q^n}{n!}=\frac{1}{1 - t + t \sqrt{1 - 2 q}},
\end{equation}
 which was previously obtained by Kuba and Panholzer \cite{Kuba-Panholzer-2007} in their study of degree distributions in increasing trees. 
 
Combining Theorem \ref{main1} with \eqref{S(t,q)-gfa} and \eqref{S(t,q)-gf} yields that   

\begin{thm}\label{main2} 
   We have 
  \begin{align}\label{gf-plane-improper}
     & \sum_{n=0}^\infty P_n(x,y) \frac{q^n}{n!} = \frac{1}{\sqrt {1-2(x+y)q}},\\[5pt]
      &  \sum_{n=0}^\infty O_n(x,y,t) \frac{q^n}{n!} = \frac{x + y}{x + y - t + t \sqrt{1 - 2 (x + y) q}}.\label{gf-plane-improper2}
  \end{align}
\end{thm}

The involution constructed in the proof of Theorem~\ref{main1} can also be used to establish symmetry properties of the following multivariate generating polynomial. For \(n\geq 1\), define
\begin{align}
    &G_n(x,y;u_1,u_2,u_3,u_4,u_5,u_6) \nonumber \\[5pt]
    &\quad =
    \sum_{T\in \mathcal{P}_n}
    x^{{\rm impr}(T)}
    y^{{\rm prop}(T)}
    u_1^{{\rm fsleaf}(T)}
    u_2^{{\rm etleaf}(T)}
    u_3^{{\rm ysleaf}(T)}
    u_4^{{\rm entleaf}(T)}
    u_5^{{\rm yerleaf}(T)}
    u_6^{{\rm ryint}(T)}. \label{defi:Gn}
\end{align} 
The additional statistics appearing in the exponents will be defined precisely in Section~\ref{sec:sym}. We obtain the following symmetry identity.

\begin{prop}\label{prop:G-sym}
For \(n\geq 2\), we have
\begin{align}
G_n(x,y;u_1,u_2,u_3,u_4,u_5,u_6)
=
G_n(y,x;u_2,u_1,u_4,u_3,u_6,u_5).
\label{sym1}
\end{align}
\end{prop}

After specializing \(x=y=1\), the polynomial \(G_n(x,y;u_1,u_2,u_3,u_4,u_5,u_6)\) satisfies a stronger family of symmetries among the refined leaf and interior-vertex statistics.

\begin{prop}\label{prop:G-sym-specialized}
For \(n\geq 2\), the polynomial 
$G_n(1,1;u_1,u_2,u_3,u_4,u_5,u_6)$ 
is invariant under the action of the subgroup 
\[
H=\langle (1\,2),(3\,4),(5\,6)\rangle \leq S_6,
\]
where \(S_6\) acts on the variables \(u_1,\ldots,u_6\) by permuting their indices; that is, for every \(\sigma\in H\),
\[
G_n(1,1;u_{\sigma(1)},u_{\sigma(2)},u_{\sigma(3)},
u_{\sigma(4)},u_{\sigma(5)},u_{\sigma(6)})
=
G_n(1,1;u_1,u_2,u_3,u_4,u_5,u_6).
\] 
\end{prop}

Recently, Dong, Du, Ji and Zhang~\cite{Dong-Du-Ji-Zhang-2025} established several symmetry identities for statistics on plane trees and tip-augmented plane trees. These identities were later proved bijectively by Li and Lin~\cite{Li-Lin-2025}. We show that these symmetries arise naturally as specializations of Propositions~\ref{prop:G-sym} and~\ref{prop:G-sym-specialized}, thereby providing a unified involutive proof of these results.

We conclude the introduction by explaining how the results of Theorem~\ref{main1} can be reformulated in terms of Stirling permutations and quasi-Stirling permutations. 
Note that  $S_{n,r}$ appearing in Theorem \ref{main1} also  counts  Stirling permutations on $\{1,1,2,2,\ldots,n,n\}$ with $r$ blocks.  For $n\geq 1$, let $[n]_2$ denote the multiset $\{1,1,2,2, \ldots, n, n\}$.   Recall that a Stirling permutation on $[n]_2$,   introduced by Gessel and Stanley \cite{Gessel-Stanley-1978}, is  a permutation $\sigma=\sigma_1\sigma_2\cdots \sigma_{2n}$ on $[n]_2$ such that,  for each $i$, the entries between the two occurrences of $i$ in $\sigma$, if any, are greater than $i$.  For example, $\sigma=6\,6\,3\,4\,5\,5\,4\,3\,1\,1\,2\,7\,7\,2$ is  a Stirling permutation on $[7]_2$.  A block of a Stirling permutation $\sigma=\sigma_1\cdots \sigma_{2n}$ on $[n]_2$ is a substring $a_i~ \cdots~ a_{j} $ with $a_i=a_j$ that is not contained in any larger such substring.  For example, the Stirling permutation $\sigma=6\,6\,3\,4\,5\,5\,4\,3\,1\,1\,2\,7\,7\,2$  has four blocks as seen in its block decomposition: $(6\,6)(3\,4\,5\,5\,4\,3)(1\,1)(2\,7\,7\,2)$. For more details, please see   Janson, Kuba and Panholzer \cite{Janson-Kuba-Panholzer-2011} and  Remmel and Wilson \cite{Remmel-Wilson-2015}.

Janson~\cite{Janson-2008} and Koganov~\cite{Koganov-1996} independently  presented a bijection between increasing plane trees with \(n\) edges and Stirling permutations of the multiset \([n]_2\). We refer to this bijection as the Koganov--Janson bijection. Under this bijection, increasing plane trees with \(n\) edges and exactly \(r\) children of the root \(0\) are mapped bijectively to Stirling permutations of \([n]_2\) with exactly \(r\) blocks.

 More precisely,  given an increasing plane tree rooted at $0$, perform a depth-first walk starting at the root and visiting children from left to right. Each time an edge is traversed, record the label of its non-root vertex incident to this edge. Since every edge is traversed once away from the root and once towards the root, each label in $[n]$ is recorded twice. The resulting word is a Stirling permutation on $[n]_2$. Moreover, the $r$ children of the root $0$ determine the $r$ distinct
blocks of the resulting permutation, thereby preserving the correspondence
between the root degree and the number of blocks.

 For the increasing plane tree  depicted in Fig. \ref{An example for increasing plane tree}, applying the Koganov--Janson correspondence yields the Stirling permutation:
 $$1\,2\,7\,7\,3\,5\,5\,4\,4\,3\,2\,1\,6\,6,$$
 which  has two blocks as shown in its block decomposition:   $$(1\,2\,7\,7\,3\,5\,5\,4\,4\,3\,2\,1)(6~6).$$

 Archer, Gregory, Pennington, and Slayden~\cite{Archer-Gregory-Pennington-Slayden-2019} introduced quasi-Stirling permutations by showing that the Koganov--Janson correspondence between Stirling permutations and increasing plane trees extends to a correspondence between quasi-Stirling permutations and labeled plane trees. As observed by Elizalde~\cite{Elizalde-2021}, these permutations can also be represented as noncrossing matchings on \([2n]\) with arcs labeled by \(1,2,\ldots,n\).

Recall that a {\it quasi-Stirling permutation} is a word
\(\sigma=\sigma_1\sigma_2\cdots\sigma_{2n}\) obtained by permuting the multiset \([n]_2\) and avoiding the patterns \(1212\) and \(2121\). Let \(\mathcal Q_n\) denote the set of quasi-Stirling permutations on \([n]_2\). Stirling permutations form a subclass of quasi-Stirling permutations: they are obtained by imposing the stronger condition that, for each \(i\), all entries between the two occurrences of \(i\) are larger than \(i\).

Let \(\sigma=\sigma_1\sigma_2\cdots\sigma_{2n}\in \mathcal Q_n\). If \(p<q\) are the positions of the two occurrences of \(i\), then the subword
\(\sigma_p\sigma_{p+1}\cdots\sigma_q\) is called the {\it \(i\)-arc} and is denoted by \(B(i)\).

For example, let
\begin{equation}\label{exam:quasi}
\pi=5\,1\,1\,2\,2\,5\,4\,8\,3\,3\,8\,6\,6\,4\,7\,7\in\mathcal Q_8.
\end{equation}
Its arc diagram is shown below:
\begin{equation*} 
\begin{tikzpicture}[
    baseline=(base),
    x=0.42cm,
    y=0.45cm,
    every node/.style={font=\small},
    arc/.style={line width=0.45pt}
]
\coordinate (base) at (8.5,0);

\foreach \i/\a in {
1/5,2/1,3/1,4/2,5/2,6/5,7/4,8/8,
9/3,10/3,11/8,12/6,13/6,14/4,15/7,16/7}
{
    \node (p\i) at (\i,0) {\(\a\)};
}

\draw[arc] (p1) to[out=90,in=90,looseness=0.75] (p6);
\draw[arc] (p2) to[out=90,in=90,looseness=1.7] (p3);
\draw[arc] (p4) to[out=90,in=90,looseness=1.7] (p5);

\draw[arc] (p7) to[out=90,in=90,looseness=0.7] (p14);
\draw[arc] (p8) to[out=90,in=90,looseness=0.95] (p11);
\draw[arc] (p9) to[out=90,in=90,looseness=1.7] (p10);
\draw[arc] (p12) to[out=90,in=90,looseness=1.7] (p13);
\draw[arc] (p15) to[out=90,in=90,looseness=1.7] (p16);

\end{tikzpicture}
\end{equation*}
From this diagram, we have \(B(i)=ii\) for \(i=1,2,3,6,7\), and
\[
B(5)=5\,1\,1\,2\,2\,5,\qquad
B(4)=4\,8\,3\,3\,8\,6\,6\,4,\qquad
B(8)=8\,3\,3\,8.
\]

A {\it block} of a quasi-Stirling permutation is an arc that is not contained in any larger arc. More generally, if \(B(i)\) is an arc, then a {\it sub-block} of \(B(i)\) is an arc \(B(j)\) that is properly contained in \(B(i)\) and is maximal with respect to this property; equivalently, there is no arc \(B(k)\) satisfying
\[
B(j)\subsetneq B(k)\subsetneq B(i).
\]
Thus \(B(i)\) admits a unique decomposition
\[
B(i)=i\,B(j_1)B(j_2)\cdots B(j_k)\,i,
\]
where \(B(j_1),\ldots,B(j_k)\) are the sub-blocks of \(B(i)\), listed from left to right. This decomposition may be applied recursively to each sub-block.

For the permutation \(\pi\) above, the block decomposition is
\[
(5\,1\,1\,2\,2\,5)\,(4\,8\,3\,3\,8\,6\,6\,4)\,(7\,7).
\]
The nontrivial sub-block decompositions are
\[
B(5)=5(11)(22)5,\qquad
B(4)=4(8338)(66)4,\qquad
B(8)=8(33)8.
\]

We now define improper and proper arcs of quasi-Stirling permutations. Recall that a right-to-left minimum of a word is an entry that is strictly smaller than every entry to its right. Let $\pi$ be a quasi-Stirling permutation in $\mathcal Q_n$ and let $B(i)$ be a block or a sub-block of  $\pi$ with nontrivial block decomposition: 
\[
B(i)=i\,B(j_1)B(j_2)\cdots B(j_k)\,i
\]
For $1\leq s\leq k$, the sub-block \(B(j_s)\)   is called {\it improper} in $B(i)$   if it contains at least one right-to-left minimum of the block \(B(i)\); otherwise it is called  {\it proper} in $B(i)$.

For the quasi-Stirling permutation $\pi$ in \eqref{exam:quasi}, the right-to-left minima 
of $B(5)$ that lie inside its sub-blocks are contained in $B(1)$ and $B(2)$.
Hence $B(1)$ and $B(2)$ are improper in $B(5)$. Similarly, $B(8)$ contains 
a right-to-left minimum of $B(4)$, while $B(6)$ does not; 
hence $B(8)$ is improper in $B(4)$, whereas $B(6)$ is proper 
in $B(4)$. Finally, since $B(3)$ contains a right-to-left minimum of $B(8)$, 
the sub-block $B(3)$ is improper in $B(8)$.

Let \({\rm impr}(\pi)\), \({\rm prop}(\pi)\), and \({\rm block}(\pi)\) denote the numbers of improper arcs, proper arcs, and blocks of \(\pi\), respectively. Thus, for the quasi-Stirling permutation \(\pi\) in \eqref{exam:quasi}, we have,
\[
{\rm impr}(\pi)=4,\qquad {\rm prop}(\pi)=1,\qquad {\rm block}(\pi)=3.
\]

For a Stirling permutation, every sub-block arc is proper. Conversely, any quasi-Stirling permutation with no improper arcs must be a Stirling permutation.

Archer et al.~\cite{Archer-Gregory-Pennington-Slayden-2019} observed that the bijection $\varphi$ also provides a bijection between  labeled plane trees rooted at \(0\) and quasi-Stirling permutations.   Starting from the root \(0\), we recursively visit the children of each vertex  from left to right. Whenever the traversal moves from a vertex to a child labeled \(i\), we record \(i\); when it returns from this child to its parent, we record \(i\) again. The resulting word is a quasi-Stirling permutation. Conversely, the noncrossing matching representation uniquely determines the labeled plane tree. For example, the quasi-Stirling permutation in  \eqref{exam:quasi} corresponds to the labeled plane tree as shown in Fig. \ref{Fig. An example for labeled plane treeII}. 

 \begin{figure}[H]
    \centering
 \begin{tikzpicture}
 [scale=0.6,vertex/.style={shape=circle, draw, inner sep=1.5pt, fill=black},
 subtree/.style={shape=ellipse, draw,minimum width=1.5cm, minimum height=.5cm},
 every fit/.style={ellipse,draw,inner sep=-2pt},
sibling distance=2.5cm,level distance=14mm,
 leaf/.style={label={[name=#1]below:$ $}},auto]
 
\node[vertex,label=0:{$0$}]{}[grow=down, sibling distance=2cm,level distance=15mm]
  child {node [vertex, label=-180:{$5$}]{} [grow=down, sibling distance=1.2cm,level distance=15mm]
  child {node [vertex, label=-90:{$1$}]{} edge from parent[solid] }
  child {node [vertex, label=-90:{$2$}]{} edge from parent[solid] }edge from parent[solid]}
  child {node [vertex, label=0:{$4$}]{}[sibling distance=12mm,level distance=15mm]
  child{node[vertex,label=0:{$8$}] {}child{node[vertex,label=0:{$3$}] {}edge from parent[solid]} edge from parent[solid] } child{node[vertex,label=0:{$6$}] {} edge from parent[solid]}}
  child{node[vertex,label=-90:{$7$}]{} edge from parent[solid] };
\end{tikzpicture}
\caption{A labeled plane tree with root $0$ and $8$ edges.}
\label{Fig. An example for labeled plane treeII}
    \end{figure}

Under the bijection, each block $B(i)$ corresponds to  the child of the root labeled $0$, and the sub-blocks of \(B(j)\) correspond to  the children of the vertex labeled \(j\), in their left-to-right order. Moreover, the smallest entry contained in a sub-block is precisely the smallest label in the corresponding subtree. Consequently, an arc is improper if and only if the corresponding edge in the labeled plane tree is improper.
 
Define 
\begin{equation}
    Q_n(x,y,t)=\sum_{\pi\in \mathcal{Q}_{n}} x^{{\rm impr}(\pi)} y^{{\rm prop}(\pi)} t^{{\rm block}(\pi)}
\end{equation}
with the convention that $Q_0(x,y,t)=1$. The above correspondence immediately yields the following proposition.
\begin{prop}\label{prop: OQ}
For $n\geq 0$, 
\begin{equation}
    Q_{n}(x,y,t)=O_{n}(x,y,t).
\end{equation}
\end{prop}
 Combining Proposition \ref{prop: OQ} with \eqref{gf-plane-improper2}, we arrive at
 \begin{equation}
     \sum_{n=0}^\infty Q_n(x,y,t) \frac{q^n}{n!} = \frac{x + y}{x + y - t + t \sqrt{1 - 2 (x + y) q}}.
 \end{equation}


The remainder of this paper is organized as follows. In Section \ref{sec:involution}, we introduce the core involution and discover its commutativity. Specifically,  the transformation is independent of the order in which improper edges are processed. Section \ref{sec:th1proof} is devoted to the proof of Theorem \ref{main1}, establishing the polynomial analogue of the identity $(n+1)!C_n = 2^n(2n-1)!!$.  
Finally, in Section~\ref{sec:sym}, we use the involution constructed in Section~\ref{sec:involution} to prove symmetry properties of the polynomial
\(G_n(x,y;u_1,u_2,u_3,u_4,u_5,u_6)\) defined in \eqref{defi:Gn} and of its specializations.

 \section{An involution on labeled plane trees} \label{sec:involution}

 For a nonempty labeled plane tree $T \in \mathcal{P}_{n}$, let $e=(i,j)$ be an edge of $T$. Assume that $i$ has $p$ children, that is, $k_1,\ldots, k_{t-1}, j, k_{t+1}, \ldots, k_p$, where $j$ is its $t$-th child from left to right. Let $\tau_i$ denote the subtree rooted by $i$.  
 The plane tree $T$ is divided into three parts according to the edge $e=(i,j)$ (see Fig. \ref{treeABC}): 
 \begin{itemize}
 \item $A_e$: the forest formed by subtrees $\tau_{k_1},\ldots, \tau_{k_{t-1}}$;
 
 \item $B_e$: the forest formed by  subtrees $\tau_{l_1},\ldots, \tau_{l_q}$;
 
 \item $C_e$: the forest formed by subtrees $\tau_{k_{t+1}},\ldots,\tau_{k_{p}}$.
 \end{itemize}
 
 Note that the forests  $A_e$, $B_e$ and $C_e$ can be empty.  Let $\hat{T}$  be the labeled plane tree obtained from $T$   by  swapping $\{{i},C_e\}$ and $\{{j},B_e\}$ (see Fig. \ref{treeACB}) and set  $\phi_e(T)=\hat{T}$. 
 
 It is clear that the map $\phi_{e}$ is an involution on $\mathcal{P}_n$. Moreover, the sets of edges in $T$  and $\hat{T}$ are  identical.  The only differences between the two trees lie in the order of the edges and the vertex connections of the edge $e$. 
\begin{figure}[H]
    \centering
    
      \begin{tikzpicture}
[vertex/.style={shape=circle, draw, inner sep=1pt, fill=black},
subtree/.style={shape=ellipse, draw,minimum width=1.5cm, minimum height=.5cm},
every fit/.style={ellipse,draw,inner sep=-2pt},
sibling distance=2cm,level distance=1cm,
leaf/.style={label={[name=#1]below:$ $}}]

\node{}[grow=down]
  child {node [vertex, label=30:{$i$}]{}[level distance=12mm]
   child {node [vertex, label=180:{$k_1$}] (a's parent) {}  [level distance=12mm]
    child {node [subtree, rotate=90,leaf=a,label={[yshift=-10mm]0:{\small$\tau_{k_1}$}}, label={[yshift=-10mm,xshift=10mm]0:{$\cdots$}}]{} [level distance=9mm] child {node [ rotate=70,leaf=c]{} edge from parent [solid]}edge from parent [solid]}edge from parent [solid]}
   child {node [vertex, label={[yshift=-12mm,xshift=-6mm]0:{\small$\tau_{k_{t-1}}$}}, label=180:{{$k_{t-1}$}}] (b's parent) {}  [level distance=12mm]
   child {node [subtree, rotate=90,leaf=b]{} edge from parent [solid]}edge from parent [solid]}   
   child {node [vertex,label=180:{$j$}]{} [level distance=12mm]
   child {node [vertex, label={[yshift=-12mm,xshift=-4mm]0:{\small$\tau_{l_{1}}$}},label=0:{$l_1$}] (z's parent) {}  [level distance=12mm]
    child {node [subtree, rotate=90,leaf=z,label={[yshift=-10mm,xshift=10mm]0:{$\cdots$}}]{} [level distance=9mm] child {node [ rotate=70,leaf=v]{} edge from parent [solid]}edge from parent [solid]}edge from parent [solid]}
   child {node [vertex, label={[yshift=-12mm,xshift=-4mm]0:{\small$\tau_{l_{q}}$}}, label=180:{$l_{q}$}] (x's parent) {}  [level distance=12mm]
   child {node [subtree, rotate=90,leaf=x]{} edge from parent [solid]}edge from parent [solid]} 
   edge from parent [solid]node[right]{$e$}} 
   child {node [vertex, label={[yshift=-12mm,xshift=-6mm]0:{\small$\tau_{k_{t+1}}$}}, label=0:{$k_{t+1}$}] (d's parent) {}  [level distance=12mm]
   child {node [subtree, rotate=90,leaf=d,label={[yshift=-10mm,xshift=10mm]0:{$\cdots$}}]{}[level distance=9mm] child {node [ rotate=70,leaf=f]{}edge from parent [solid] }edge from parent [solid] }edge from parent [solid]}
   child {node [vertex, label={[yshift=-12mm,xshift=-4.5mm]0:{\small$\tau_{k_{p}}$}},label=0:{$k_p$}] (e's parent) {}  [level distance=12mm]
   child {node [subtree, rotate=90,leaf=e]{} edge from parent [solid]}edge from parent [solid] }
   edge from parent [dashed]};
    \node[draw,dashed,xshift=-2mm,fit=(a's parent)(a) (b) (c)  ,label=below:$A_e$] {}; 
     \node[draw,dashed,xshift=-2mm,fit=(z's parent)(z) (x) (v)  ,label=below:$B_e$] {}; 
      \node [dashed,xshift=-2mm,fit=(e) (f) (d's parent),label=below:$C_e$] {};
\end{tikzpicture} 
        \caption{ The decomposition of a labeled plane tree  $T$ based on edge  $e$}
         \label{treeABC}
\end{figure}
  
\begin{figure}[H]
    \centering
\caption*{\hspace{-8cm}\vspace{-1cm}$T$} 
\caption*{\hspace{8cm} $\phi_e(T)$} 
\begin{tikzpicture}
[vertex/.style={shape=circle, draw, inner sep=1pt, fill=black},
subtree/.style={shape=ellipse, dashed, draw,minimum width=1.7cm, minimum height=.8cm},
every fit/.style={ellipse,draw,inner sep=-2pt},
sibling distance=1.8cm,level distance=1cm,
leaf/.style={label={[name=#1]below:$ $}}]

\node{}[grow=down]
  child {node [vertex,name=1,label=30:{$i$}]{}[level distance=12mm]
    child {node [subtree, rotate=30,leaf=a,label={[yshift=-7mm,xshift=-12mm]0:{\small$A$}}]{}  edge from parent [solid]}   
   child {node [vertex,name=2, label={[yshift=-12mm,xshift=-4mm]0:{\small$B$}}, label=180:{$j$}]  {}  [level distance=12mm]
   child {node [subtree, rotate=90,leaf=x]{} edge from parent [solid]}
   edge from parent [solid]node[right]{$e$}} 
   child {node [subtree, rotate=-30,leaf=a,label={[yshift=8mm,xshift=-12mm]0:{\small$C$}}]{} [level distance=9mm] edge from parent [solid]}  
   edge from parent [dashed]};
        \draw[->, line width=1pt, shorten >=2pt, shorten <=2pt] (3,-2) -- (5,-2) 
              node[yshift=-7.5mm,xshift=-10mm] {$\phi_{\hat{e}}$};
              \draw[->, line width=1pt, shorten >=2pt, shorten <=2pt] (5,-2.3) -- (3,-2.3) 
              node[yshift=7.5mm,xshift=10mm] {$\phi_e$};
              
\end{tikzpicture}
\begin{tikzpicture}
[vertex/.style={shape=circle, draw, inner sep=1pt, fill=black},
subtree/.style={shape=ellipse, dashed, draw,minimum width=1.7cm, minimum height=.8cm},
every fit/.style={ellipse,draw,inner sep=-2pt},
sibling distance=1.8cm,level distance=1cm,
leaf/.style={label={[name=#1]below:$ $}}]

\node{}[grow=down]
  child {node [vertex, label=30:{$j$}]{}[level distance=12mm]
    child {node [subtree, rotate=30,leaf=a,label={[yshift=-7mm,xshift=-12mm]0:{\small$A$}}]{}  edge from parent [solid]}   
   child {node [vertex, label={[yshift=-12mm,xshift=-4mm]0:{\small$C$}}, label=180:{$i$}]  {}  [level distance=12mm]
   child {node [subtree, rotate=90,leaf=x]{} edge from parent [solid]}
   edge from parent [solid]node[right]{$\hat{e}$}} 
   child {node [subtree, rotate=-30,leaf=a,label={[yshift=8mm,xshift=-12mm]0:{\small$B$}}]{} [level distance=9mm] edge from parent [solid]}  
   edge from parent [dashed]};
\end{tikzpicture}
\caption{A visual illustration of the involution  $\phi$}
        \label{treeACB}
\end{figure}

 The main property of this involution $\phi$ is contained in the following  proposition:

\begin{prop}\label{imppro}
For a labeled plane tree $T$ and an edge  $e$ in $T$, define $\hat{T}:=\phi_{e}(T)$. The proper or improper status of the edge $e$  is reversed between $T$ and $\hat{T}$, while all edges other than $e$  preserve their proper or improper property.
\end{prop}
\begin{proof}

Let $a=(a_1,a_2)$ be any edge of $T$.  As illustrated in Fig. \ref{treeABC},  we decomposed the plane tree $T$ into three parts based on the edge $a=(a_1,a_2)$.  Let $A_T(a), B_T(a), C_T(a)$ be the forests formed by subtrees of $T$, as defined in the construction of the involution $\phi$. Define 
 $\mathcal{B}_T(a)=\{a_2,B_T(a)\}$ and $\mathcal{C}_T(a)= \{a_1,C_T(a)\}$.  By definition,  if $a$ is a proper edge of $T$, then $\min \mathcal{B}_T(a) > \min \mathcal{C}_T(a),$  otherwise, $\min \mathcal{B}_T(a) < \min \mathcal{C}_T(a).$ 
 
 Similarly, for $\hat{T}$,  we define $\mathcal{B}_{\hat{T}}(a)$ and $\mathcal{C}_{\hat{T}}(a)$ with respect to  the edge $a$ in $\hat{T}$. The proper/improper status of $a$ in $\hat{T}$ is also determined by comparing 
 the minima of $\mathcal{B}_{\hat{T}}(a)$ and $\mathcal{C}_{\hat{T}}(a)$. 

 From the construction of the involution $\phi_e$, we see 
 
 if $a=e$, then $\mathcal{B}_T(e)=\mathcal{C}_{\hat{T}}({e})$ and $\mathcal{C}_T(e)=\mathcal{B}_{\hat{T}}({e})$. This implies the proper/improper status of  $e$  is reversed between $T$ and $\hat{T}$. 

if $a\neq e$, then $\mathcal{B}_T(a)=\mathcal{B}_{\hat{T}}({a})$ and $\mathcal{C}_T(a)=\mathcal{C}_{\hat{T}}({a})$ Thus, the proper/improper status of the edge $a$  remains unchanged in $T$ and $\hat{T}$.    This completes the proof.   
\end{proof}

From the construction of $\phi$, it is not hard to check that 
\begin{prop}\label{commutativity} Given a labeled plane tree $T$,  the involutions $\phi$ commute for  any pair of edges $e_1,e_2$ in $T$: 
\[\phi_{e_2}\phi_{e_1}(T)=\phi_{e_1}\phi_{e_2}(T).\]
\end{prop}

 \section{Proof of Theorem \ref{main1}} \label{sec:th1proof}

 In this section, we give a proof of Theorem \ref{main1} based on  the involution $\phi$ defined in Section 2.

\proof[Proof of Theorem \ref{main1}] Let $\mathcal{I}_n$ denote the set of increasing plane trees with $n$ edges. To prove \eqref{Pnxy},  we first define edge labelings for  plane trees in  $\mathcal{P}_n$ and  $\mathcal{I}_n$. 

Let $T$ be a labeled plane tree with $n$ edges. We assign a label to each edge of   $T$: improper edges are labeled by $x$, and proper edges are labeled by $y$. The weight of $T$, denoted by ${\rm wt}(T)$ is defined to be the product of these edge labels. It is clear  that,  
\[P_n(x,y)=\sum_{T \in \mathcal{P}_n} {\rm wt}(T).\]
For the plane tree $T\in \mathcal{I}_n$, we introduce a free labeling:   each edge may be labeled either $x$ or $y$. The weight of $T \in \mathcal{I}_n$, denoted by $\widetilde{{\rm wt}}(T)$, is the product of its edge labels. It follows that 
\[\sum_{T \in \mathcal{I}_n}\widetilde{{\rm wt}}(T)=(2n-1)!!(x+y)^n.\]
Proving \eqref{Pnxy} is  thus equivalent to showing
\begin{equation}\label{pfaa}
\sum_{T \in \mathcal{P}_n} {\rm wt}(T)=\sum_{T \in \mathcal{I}_n}\widetilde{{\rm wt}}(T).
\end{equation}
Let $T$ be a labeled plane tree in $\mathcal{P}_n$ with $k$ improper edges. These improper edges are denoted by $e_1, e_2, \ldots, e_k$  in the order they appear in the depth-first walk of the rooted plane tree. Applying the involutions  $\phi_{e_1}, \ldots, \phi_{e_k}$  to $T$ successively  (preserving edge labels), we get the labeled plane tree $\hat{T}=\phi_{e_k}\phi_{e_{k-1}}\cdots \phi_{e_1}(T)$. Define $\Phi=\phi_{e_k}\phi_{e_{k-1}}\cdots \phi_{e_1}$, so that   $\hat{T}=\Phi(T)$, see Fig. \ref{fig:phi}. 
By Proposition \ref{imppro},  $\hat{T}$ contains no improper edges, and so $\hat{T} \in \mathcal{I}_n$ (the set of increasing plane trees with $n$ edges). Moreover,  ${\rm wt}(T)=\widetilde{{\rm wt}}(\hat{T})$, where the $x$-labeled edges in $\widetilde{{\rm wt}}(\hat{T})$ correspond to  the improper edges of $T$. 

Conversely, let $\hat{T} \in \mathcal{I}_n$ with the weight $\widetilde{{\rm wt}}(\hat{T})$.  Assume that there are $k$ edges labeled by $x$. Let  $\hat{e}_1, \hat{e}_2, \ldots \hat{e}_k$ be its  $x$-labeled edges in the depth-first walk order. By successively applying the involutions $\phi_{\hat{e}_1}, \ldots, \phi_{\hat{e}_k}$ to $\hat{T}$   (preserving labels), we  obtain the labeled plane tree ${T}=\phi_{\hat{e}_k}\phi_{\hat{e}_{k-1}}\cdots \phi_{\hat{e}_1}(\hat{T})$. Define $\Psi=\phi_{\hat{e}_k}\phi_{\hat{e}_{k-1}}\cdots \phi_{\hat{e}_1}$, which gives ${T}=\Psi(\hat{T}).$ 
By Proposition \ref{imppro}, ${T}$ has exactly   $k$ improper edges (all labeled $x$), and so ${\rm wt}(T)=\widetilde{{\rm wt}}(\hat{T})$. Furthermore, by Proposition  \ref{commutativity},  we derive that  $T=\Psi(\Phi(T))$ for all $T \in \mathcal{P}_n$ and $\hat{T}=\Phi(\Psi(\hat{T}))$ for all $\hat{T} \in \mathcal{I}_n$. This establishes a weight-preserving bijection between the two sums in \eqref{pfaa}, completing the proof of  \eqref{pfaa}   and thus of  \eqref{Pnxy}.

 The proof of \eqref{Fnxyt}  is similar to that of  \eqref{Pnxy}. It suffices to modify the edge labelings for plane trees in  $\mathcal{O}_n$ and  $\mathcal{I}_n$ as follows:  For any $T \in \mathcal{O}_n$,  label the edges  connected to vertex 0 with $t$ and  label all other edges with 
  $x$ or $y$ using the same rule as in 
 $\mathcal{P}_n$ (i.e., $x$  for improper edges, $y$ for proper edges). The weight of $T \in \mathcal{O}_n$, denoted by ${\rm wt}^{r}(T)$, is defined to be the product of these edge labels. Clearly,  
\[O_n(x,y, t)=\sum_{T\in \mathcal{O}_{n}}{\rm wt}^r(T).\]
For increasing plane trees $T\in \mathcal{I}_n$,  we introduce a modified free labeling: Edges connected to vertex $0$ are labeled by $t$,  while all other edges may be labeled either   $x$ or $y$.  The weight of $T \in \mathcal{I}_n$, denoted by $\widetilde{{\rm wt}}^r(T)$, is the product of its edge labels. It follows that 
\[\sum_{r=1}^n t^r S_{n,r} (x+y)^{n-r}=\sum_{T \in \mathcal{I}_n}\widetilde{{\rm wt}}^r(T).\]
Thus, establishing \eqref{Fnxyt}  reduces to showing that
\begin{equation}\label{pfbb}
\sum_{T\in \mathcal{O}_{n}}{\rm wt}^r(T)=\sum_{T \in \mathcal{I}_n}\widetilde{{\rm wt}}^r(T).
\end{equation}
The bijection $\Phi$ used in the proof of  \eqref{pfaa} can be readily adapted to establish \eqref{pfbb}, and consequently   \eqref{Fnxyt}.  
 This completes the proof of Theorem \ref{main1}.  \qed 

 \begin{figure}[H]
 \begin{tikzpicture}
 [scale=0.7,vertex/.style={shape=circle, draw, inner sep=1.5pt, fill=black},
 subtree/.style={shape=ellipse, draw,minimum width=1.5cm, minimum height=.5cm},
 every fit/.style={ellipse,draw,inner sep=-2pt},
sibling distance=3.5cm,level distance=15mm,
 leaf/.style={label={[name=#1]below:$ $}},auto]
 
\node[vertex,name=5,label=0:{$1$}]{}[grow=down]
  child {node [vertex, name=1, label=-180:{$6$}]{} 
  child {node [vertex, label=-90:{$0$}]{} 
  edge from parent[solid]node[xshift=-6mm]{$e_2$}node[right]{\footnotesize$(x)$} }
  edge from parent[solid]node[yshift=4mm,xshift=-6mm]{$e_1$} node[yshift=3mm,xshift=-1mm]{\footnotesize$(x)$}}
  child {node [vertex, label=0:{$3$}]{}[sibling distance=15mm,level distance=15mm]
  child{node[vertex,label=-90:{$7$}] {}edge from parent[solid]  node[left]{\footnotesize$(y)$}} child{node[vertex,label=-90:{$2$}] {}edge from parent[solid]node[yshift=1mm,xshift=-2mm]{\footnotesize$(x)$} node[yshift=1mm,xshift=-6.5mm]{$e_3$}} child{node[vertex,label=-90:{$5$}] {}edge from parent[solid]node[yshift=-5mm,xshift=0mm]{\footnotesize$(y)$}} child{node[vertex,label=-90:{$4$}]{} edge from parent[solid]node[right]{\footnotesize$(y)$} }edge from parent[solid]node[right]{\footnotesize$(y)$}};
  \draw[red, thick] (5) -- (1);
  \draw[->, line width=1pt, shorten >=2pt, shorten <=2pt]  (5,-1.3) --(8,-1.3)
    node[yshift=4mm,xshift=-10mm] {$\phi_{e_1}$};
\end{tikzpicture}
    \begin{tikzpicture}
 [scale=0.7,vertex/.style={shape=circle, draw, inner sep=1.5pt, fill=black},
 subtree/.style={shape=ellipse, draw,minimum width=1.5cm, minimum height=.5cm},
 every fit/.style={ellipse,draw,inner sep=-2pt},
sibling distance=3cm,level distance=15mm,
 leaf/.style={label={[name=#1]below:$ $}},auto]
 
\node[vertex,name=1, label=0:{$6$}]{}[grow=down]
  child {node [vertex, name=5,  label=-180:{$1$}]{}  
  child {node [vertex, name=3, label=0:{$3$}]{}[sibling distance=18mm,level distance=15mm] 
  child {node [vertex, label=-90:{$7$}]{}edge from parent[solid]node[yshift=3mm,xshift=-8mm]{\footnotesize$(y)$}} 
  child{node[vertex,label=-90:{$2$}] {}edge from parent[solid]node[yshift=2mm,xshift=-2mm]{\footnotesize$(x)$}node[yshift=2mm,xshift=-6mm]{$e_3$}} 
  child{node[vertex,label=-90:{$5$}] {}edge from parent[solid]node[yshift=-4mm,xshift=-1mm]{\footnotesize$(y)$}} 
  child{node[vertex,label=-90:{$4$}] {} edge from parent[solid]node[yshift=-4mm,xshift=2mm]{\footnotesize$(y)$} }edge from parent[solid]node[right]{\footnotesize$(y)$}}edge from parent[solid]node[left]{\footnotesize$(x)$}}
  child{node[vertex,name=7, label=-90:{$0$}]{}edge from parent[solid]node[right]{\footnotesize$(x)$}node[left]{$e_2$}};
  
  \draw[red, thick] (1) -- (7);
   \draw[->, line width=1pt, shorten >=2pt, shorten <=2pt]  (2,-2.7) --(5,-2.7)
        node[yshift=4mm,xshift=-10mm] {$\phi_{e_2}$};
\end{tikzpicture}
\begin{tikzpicture}
 [scale=0.7,vertex/.style={shape=circle, draw, inner sep=1.5pt, fill=black},
 subtree/.style={shape=ellipse, draw,minimum width=1.5cm, minimum height=.5cm},
 every fit/.style={ellipse,draw,inner sep=-2pt},
sibling distance=3cm,level distance=15mm,
 leaf/.style={label={[name=#1]below:$ $}},auto]
 
\node[vertex,name=1, label=0:{$0$}]{}[grow=down]
  child {node [vertex, name=5,  label=-180:{$1$}]{}  
  child {node [vertex, name=3, label=0:{$3$}]{}[sibling distance=18mm,level distance=15mm] 
  child {node [vertex, label=-90:{$7$}]{}edge from parent[solid]node[yshift=3mm,xshift=-8mm]{\footnotesize$(y)$}} 
  child{node[vertex, name=2, label=-90:{$2$}] {}edge from parent[solid]node[yshift=2mm,xshift=-2mm]{\footnotesize$(x)$}node[yshift=2mm,xshift=-6mm]{$e_3$}} 
  child{node[vertex,label=-90:{$5$}] {}edge from parent[solid]node[yshift=-4mm,xshift=-1mm]{\footnotesize$(y)$}} 
  child{node[vertex,label=-90:{$4$}] {} edge from parent[solid]node[yshift=-4mm,xshift=2mm]{\footnotesize$(y)$} }edge from parent[solid]node[right]{\footnotesize$(y)$} }edge from parent[solid]node[left]{\footnotesize$(x)$}}
  child{node[vertex,name=7, label=-90:{$6$}]{}edge from parent[solid]node[right]{\footnotesize$(x)$}};
  
  \draw[red, thick] (2) -- (3);
   \draw[->, line width=1pt, shorten >=2pt, shorten <=2pt]  (2.5,-2.7) --(5.5,-2.7)
         node[yshift=4mm,xshift=-10mm] {$\phi_{e_3}$};
\end{tikzpicture}\qquad 
\begin{tikzpicture}
 [scale=0.6,vertex/.style={shape=circle, draw, inner sep=1.5pt, fill=black},
 subtree/.style={shape=ellipse, draw,minimum width=1.5cm, minimum height=.5cm},
 every fit/.style={ellipse,draw,inner sep=-2pt},
sibling distance=3cm,level distance=15mm,
 leaf/.style={label={[name=#1]below:$ $}},auto]
\node[vertex,name=1,label=0:{$0$}]{}[grow=down]
  child {node [vertex, name=3, label=180:{$1$}]{}  
  child {node [vertex, name=5,label=0:{$2$}]{}
  child {node [vertex, label=-90:{$7$}]{}node[yshift=6mm,xshift=3mm]{\footnotesize$(y)$}}
  child {node [vertex, label=0:{$3$}]{}
  child {node [vertex, label=-90:{$5$}]{}node[yshift=6mm,xshift=3mm]{\footnotesize$(y)$}}
  child {node [vertex, label=-90:{$4$}]{}
  node[yshift=6mm,xshift=-2mm]{\footnotesize$(y)$}
  }node[yshift=6mm,xshift=-1mm]{\footnotesize$(x)$}}node[yshift=6mm,xshift=3mm]{\footnotesize$(y)$}} node[yshift=6mm,xshift=1mm]{\footnotesize$(x)$}}
  child {node [vertex, label=-90:{$6$}]{}edge from parent[solid]node[yshift=-3mm,xshift=1mm]{\footnotesize$(x)$}};
\end{tikzpicture} 
\caption{An example of the bijection  $\Phi=\phi_{e_1}\phi_{e_2}\phi_{e_3}$ } \label{fig:phi}
\end{figure}

\section{Proofs of Proposition~\ref{prop:G-sym} and Proposition \ref{prop:G-sym-specialized} }\label{sec:sym}

This section aims to prove  Proposition~\ref{prop:G-sym} and Proposition~\ref{prop:G-sym-specialized} with the aid of the involution in   Section~\ref{sec:involution}.

 We first recall some statistics for plane trees.  Chen, Deutsch, and Elizalde \cite{Chen-Deutsch-Elizalde-2006} classified the leaves of a plane tree into old and young leaves. Dong et al. \cite{Dong-Du-Ji-Zhang-2025} called a leaf without siblings a \textit{singleton leaf}. In our convention, a leaf with siblings is called an \textit{elder leaf} if it is the rightmost child of its parent, and a \textit{young leaf} otherwise. This convention differs from that of Dong et al. \cite{Dong-Du-Ji-Zhang-2025}, where elder leaves are defined to be leftmost children; the two conventions are equivalent under reflection of plane trees. An interior vertex is called a \textit{young interior vertex} if it is not the parent of a singleton leaf or an elder leaf.
 
We next introduce finer classifications of  interior vertices and leaves.

\begin{defi}
Let \(T\) be a labeled plane tree.

\begin{itemize}
    \item A young interior vertex \(i\) is called a \textit{singleton interior vertex} if its rightmost child \(j\) has exactly one child \(k\), and  \(k\) is a leaf.
    \item A young interior vertex that is not a singleton interior vertex is called a \textit{non-singleton young interior vertex}.
    \end{itemize}
    \end{defi}
    
    \begin{defi}
    Let $T$ be a labeled plane tree.
        \begin{itemize}
    \item A leaf  is called a \textit{first singleton leaf} if it is the unique child of the rightmost child of a singleton interior vertex.

    \item A singleton leaf that is not a first singleton leaf is called a
    \textit{young singleton leaf}.

    \item An elder leaf is called an \textit{elder twin leaf} if its immediate
    left sibling is also a leaf.

    \item An elder leaf is called an \textit{elder non-twin leaf} if its
    immediate left sibling is not a leaf.

    \item A young leaf is called a \textit{second leaf} if it is the immediate left sibling of an elder twin leaf.

    \item A young leaf that is not a second leaf is called a \textit{younger leaf}.
\end{itemize}
\end{defi}

We summarize the corresponding vertex-type statistics of a plane tree $T$ in Table~\ref{table:statistics}.
 \begin{table}[H]
\centering
\caption{Vertex-type statistics for a plane tree \(T\).}
\label{table:statistics}
\begingroup
\renewcommand{\arraystretch}{0.9}
\begin{tabular}{ll}
\toprule
\({\rm sleaf}(T)\)  & \(\#\{\text{singleton leaves in }T\}\) \\
\({\rm eleaf}(T)\)  & \(\#\{\text{elder leaves in }T\}\) \\
\({\rm yleaf}(T)\)  & \(\#\{\text{young leaves in }T\}\) \\
\({\rm fsleaf}(T)\) & \(\#\{\text{first singleton leaves in }T\}\) \\
\({\rm ysleaf}(T)\) & \(\#\{\text{young singleton leaves in }T\}\) \\
\({\rm etleaf}(T)\) & \(\#\{\text{elder twin leaves in }T\}\) \\
\({\rm entleaf}(T)\)& \(\#\{\text{elder non-twin leaves in }T\}\) \\
\({\rm secleaf}(T)\) & \(\#\{\text{second leaves in }T\}\) \\
\({\rm yerleaf}(T)\)& \(\#\{\text{younger leaves in }T\}\) \\
\({\rm sint}(T)\)   & \(\#\{\text{parents of singleton leaves in }T\}\) \\
\({\rm eint}(T)\)   & \(\#\{\text{parents of elder leaves in }T\}\) \\
\({\rm yint}(T)\)   & \(\#\{\text{young interior vertices in }T\}\) \\
\({\rm fsint}(T)\)  & \(\#\{\text{parents of first singleton leaves in }T\}\) \\
\({\rm ysint}(T)\)  & \(\#\{\text{parents of young singleton leaves in }T\}\) \\
\({\rm etint}(T)\)  & \(\#\{\text{parents of elder twin leaves in }T\}\) \\
\({\rm entint}(T)\) & \(\#\{\text{parents of elder non-twin leaves in }T\}\) \\
\({\rm gsint}(T)\)  & \(\#\{\text{singleton interior vertices in }T\}\) \\
\({\rm ryint}(T)\)  & \(\#\{\text{non-singleton young interior vertices in }T\}\) \\
\bottomrule
\end{tabular}
\endgroup
\end{table}

By definition,
\[
{\rm sleaf}(T)={\rm sint}(T),\qquad
{\rm eleaf}(T)={\rm eint}(T),
\]
and
\[
{\rm fsleaf}(T)={\rm fsint}(T)={\rm gsint}(T),\qquad
{\rm ysleaf}(T)={\rm ysint}(T).
\]
Similarly,
\[
{\rm etleaf}(T)={\rm etint}(T)={\rm secleaf}(T),\qquad
{\rm entleaf}(T)={\rm entint}(T).
\]
Consequently, for every \(T\in\mathcal{P}_n\), we have
\[
n+1 =
3{\rm fsleaf}(T)+2{\rm ysleaf}(T)
+3{\rm etleaf}(T)+2{\rm entleaf}(T)
+{\rm yerleaf}(T)+{\rm ryint}(T).
\]

  We now encode these statistics in the following multivariate generating polynomial:
\begin{align*}
    &G_n(x,y;u_1,u_2,u_3,u_4,u_5,u_6)\\[5pt]
    &=\sum_{T\in \mathcal{P}_n}x^{{\rm impr}(T)}y^{{\rm prop}(T)} u_1^{{\rm fsleaf}(T)}u_2^{{\rm etleaf}(T)}u_3^{{\rm ysleaf}(T)}u_4^{{\rm entleaf}(T)}u_5^{{\rm yerleaf}(T)}u_6^{{\rm ryint}(T)}.
\end{align*}
For example, Fig.~\ref{fig:refined-statistics} illustrates the six  vertex-type statistics appearing in \(G_n\). Each marked vertex is labeled by the variable to which it contributes. 
\begin{figure}[H]
    \centering
 \begin{tikzpicture}
 [scale=0.8,vertex/.style={shape=circle, draw, inner sep=1.5pt, fill=black},
 subtree/.style={shape=ellipse, draw,minimum width=1.5cm, minimum height=.5cm},
 every fit/.style={ellipse,draw,inner sep=-2pt},
sibling distance=1cm,level distance=14mm,
 leaf/.style={label={[name=#1]below:$0$}},auto]
 
\node[vertex,label=30:{$0$}]{}[grow=down, sibling distance=3.5cm,level distance=12mm]
  child {node [vertex, label=-180:{$1$}]{}[sibling distance=1.3cm,level distance=12mm]
  child {node [vertex, label=180:{$2$}]{}
  child {node [vertex, label=-90:{$3(u_3)$}]{}
  edge from parent[solid] }edge from parent[solid] }
 child {node [vertex, label=-180:{$4$}]{}
 child {node [vertex, label=-90:{$5(u_1)$}]{}
  edge from parent[solid] }
 edge from parent[solid]}}
  child {node [vertex, label=0:{$6$}]{}[sibling distance=12mm,level distance=12mm]
  child{node[vertex,label=90:{$7(u_6)$}] {}
  child{node[vertex,label=0:{$8$}] {} 
  child{node[vertex,label=-90:{$9$}] {} 
  edge from parent[solid] }
  child{node[vertex,label=-90:{$10(u_2)$}] {} 
  edge from parent[solid] }
  edge from parent[solid] } } child{node[vertex,label=-90:{$11(u_5)$}] {}edge from parent[solid]} child{node[vertex,label=-90:{$12$}] {}edge from parent[solid]} child{node[vertex,label=-90:{$13(u_2)$}]{} edge from parent[solid] }edge from parent[solid]}
  child{node[vertex,label=-90:{$14(u_4)$}] {} 
  edge from parent[solid] };
\end{tikzpicture}
\caption{A labeled plane tree with $14$ edges.}
\label{fig:refined-statistics}
    \end{figure}

 The following proposition records how the   involution defined in Section~\ref{sec:involution} affects  the vertex-type statistics  defined above. This proposition plays a crucial role in the proofs of Proposition~\ref{prop:G-sym} and Proposition~\ref{prop:G-sym-specialized}. 

\begin{prop} \label{setsort}  For a nonempty labeled plane tree $T \in \mathcal{P}_{n}$, let $e=(i,j)$ be an edge of $T$.    
 The plane tree $T$ is divided into three parts according to the edge $e=(i,j)$ (see Fig. \ref{treeABC}):  
\begin{enumerate}
    \item[Case (1)]  If one of $B_{e}$ and $C_e$ is empty and the other is a one-vertex set, then $\phi$ interchanges an elder twin leaf and a first singleton leaf, and also interchanges a second leaf and a singleton interior vertex. 
    \item[Case (2)] If  one of $B_e$ and $C_e$  is empty and the other has more than one vertex,  then $\phi$ interchanges a younger leaf and  a non-singleton young interior vertex.
    \item[Case (3)] If one of $B_e$ and $C_e$ has a one-vertex set and the other  has more than one vertex, then $\phi$ interchanges a young singleton leaf and an elder non-twin leaf.
     \item[Case (4)]  If $B_e$ and $C_e$ are of the same type, that is, they are both empty, both one-vertex sets or both more than one vertex, then $\phi$ preserves all the relevant statistics.
\end{enumerate}
\end{prop}

We first prove Proposition~\ref{prop:G-sym} using Proposition \ref{setsort}.

\begin{proof}[Proof of Proposition~\ref{prop:G-sym}]
Let $T\in\mathcal{P}_n$, and let $E(T)$ be its edge set. Define $$\Theta=\left(\prod_{e\in E(T)}\phi_{e}\right) (T).$$
By Proposition~\ref{commutativity}, the product is independent of the order of
the factors. Moreover, by Proposition~\ref{imppro} and
Proposition~\ref{setsort}, the map \(\Theta\) interchanges the pairs of
statistics
\[
({\rm impr},{\rm prop}),\quad
({\rm fsleaf},{\rm etleaf}),\quad
({\rm ysleaf},{\rm entleaf}),\quad
({\rm yerleaf},{\rm ryint}).
\]
This gives \eqref{sym1}. \end{proof}

Figure~\ref{fig: prop1.3} illustrates an example of the bijection \(\Theta\) used in the proof of Proposition~\ref{prop:G-sym}. In this example, the corresponding weight is transformed from \(x^3y^4u_2u_3u_5^2u_6\) to \(x^4y^3u_1u_4u_5u_6^2\).

 \begin{figure}[h]
 \centering
 \begin{tikzpicture}
 [scale=0.7,vertex/.style={shape=circle, draw, inner sep=1.5pt, fill=black},
 subtree/.style={shape=ellipse, draw,minimum width=1.5cm, minimum height=.5cm},
 every fit/.style={ellipse,draw,inner sep=-2pt},
sibling distance=4cm,level distance=15mm,
 leaf/.style={label={[name=#1]below:$ $}},auto]
 
\node[vertex,name=5,label=90:{$1(u_6)$}]{}[grow=down]
  child {node [vertex, name=1, label=-180:{$6$}]{} 
  child {node [vertex, label=-180:{$0(u_3)$}]{} 
  edge from parent[solid]node[xshift=-6mm]{$e_2 $}node[right]{\footnotesize$(x)$} }
  edge from parent[solid]node[yshift=4mm,xshift=-6mm]{$e_1$} node[yshift=3mm,xshift=-1mm]{\footnotesize$(x)$}}
  child {node [vertex, label=0:{$3$}]{}[sibling distance=20mm,level distance=15mm]
  child{node[vertex,label=-90:{$7(u_5)$}] {}edge from parent[solid] node[xshift=-9mm,yshift=2mm]{$e_3 $} node[xshift=-5mm,yshift=0mm]{\footnotesize$(y)$}} child{node[vertex,label=-90:{$2(u_5)$}] {}edge from parent[solid] node[yshift=0mm,xshift=-2mm]{\footnotesize$(x)$} node[yshift=2mm,xshift=-7mm]{$e_4$}} child{node[vertex,label=-90:{$5$}] {}edge from parent[solid] node[xshift=-5.5mm,yshift=-3mm]{$e_5 $} node[yshift=-5mm,xshift=0mm]{\footnotesize$(y)$}} child{node[vertex,label=-90:{$4(u_2)$}]{} edge from parent[solid]node[xshift=-2mm,yshift=-1mm]{$e_6 $}
  node[xshift=4mm,yshift=-4mm]{\footnotesize$(y)$} }edge from parent[solid]node[xshift=-7mm,yshift=-3mm]{$e_7 $}
  node[right]{\footnotesize$(y)$}};
  \draw[red, thick] (5) -- (1);
  \draw[->, line width=1pt, shorten >=2pt, shorten <=2pt]  (4.5,-1.3) --(7,-1.3)
    node[yshift=4mm,xshift=-10mm] {$\phi_{e_1}$};
\end{tikzpicture}
    \begin{tikzpicture}
 [scale=0.6,vertex/.style={shape=circle, draw, inner sep=1.5pt, fill=black},
 subtree/.style={shape=ellipse, draw,minimum width=1.5cm, minimum height=.5cm},
 every fit/.style={ellipse,draw,inner sep=-2pt},
sibling distance=2cm,level distance=15mm,
 leaf/.style={label={[name=#1]below:$ $}},auto]
 
\node[vertex,name=1, label=0:{$6$}]{}[grow=down]
  child {node [vertex, name=5,  label=-180:{$1$}]{}  
  child {node [vertex, name=3, label=0:{$3$}]{}[sibling distance=15mm,level distance=15mm] 
  child {node [vertex, label=-90:{$7$}]{}edge from parent[solid]  node[xshift=-9mm,yshift=2mm]{$e_3 $} } 
  child{node[vertex,label=-90:{$2$}] {}edge from parent[solid]node[yshift=2mm,xshift=-6mm]{$e_4$}} 
  child{node[vertex,label=-90:{$5$}] {}edge from parent[solid]node[xshift=-5mm,yshift=-3mm]{$e_5 $}} 
  child{node[vertex,label=-90:{$4$}] {} edge from parent[solid] node[xshift=2mm,yshift=-3mm]{$e_6 $}}edge from parent[solid] node[xshift=-5mm,yshift=0mm]{$e_7 $}}edge from parent[solid]}
  child{node[vertex,name=7, label=-90:{$0$}]{}edge from parent[solid]node[xshift=0mm,yshift=-2mm]{$e_2 $}};
  
   \draw[red, thick] (1) -- (7);
   \draw[->, line width=1pt, shorten >=2pt, shorten <=2pt]  (2.5,-2.2) --(5.5,-2.2)
        node[yshift=4mm,xshift=-10mm] {$\phi_{e_2}$};
\end{tikzpicture}\quad
\begin{tikzpicture}
 [scale=0.6,vertex/.style={shape=circle, draw, inner sep=1.5pt, fill=black},
 subtree/.style={shape=ellipse, draw,minimum width=1.5cm, minimum height=.5cm},
 every fit/.style={ellipse,draw,inner sep=-2pt},
sibling distance=2cm,level distance=15mm,
 leaf/.style={label={[name=#1]below:$ $}},auto]
 
\node[vertex, label=0:{$0$}]{}[grow=down]
  child {node [vertex,   label=-180:{$1$}]{}  
  child {node [vertex, name=3,  label=0:{$3$}]{}[sibling distance=15mm,level distance=15mm] 
  child {node [vertex, name=7,label=-90:{$7$}]{}edge from parent[solid]  node[xshift=-9mm,yshift=2mm]{$e_3 $} } 
  child{node[vertex,label=-90:{$2$}] {}edge from parent[solid]node[yshift=2mm,xshift=-6mm]{$e_4$}} 
  child{node[vertex,label=-90:{$5$}] {}edge from parent[solid]node[xshift=-5mm,yshift=-3mm]{$e_5 $}} 
  child{node[vertex,label=-90:{$4$}] {} edge from parent[solid] node[xshift=2mm,yshift=-3mm]{$e_6 $}}edge from parent[solid] node[xshift=-5mm,yshift=0mm]{$e_7 $}}edge from parent[solid]}
  child{node[vertex,  label=-90:{$6$}]{}edge from parent[solid]};

   \draw[red, thick] (7) -- (3);
   \draw[->, line width=1pt, shorten >=2pt, shorten <=2pt]  (2.5,-2.7) --(5.5,-2.7)
         node[yshift=4mm,xshift=-10mm] {$\phi_{e_3}$};
\end{tikzpicture}\quad 
\begin{tikzpicture}
 [scale=0.5,vertex/.style={shape=circle, draw, inner sep=1.5pt, fill=black},
 subtree/.style={shape=ellipse, draw,minimum width=1.5cm, minimum height=.5cm},
 every fit/.style={ellipse,draw,inner sep=-2pt},
sibling distance=2cm,level distance=15mm,
 leaf/.style={label={[name=#1]below:$ $}},auto]
\node[vertex,name=1,label=0:{$0$}]{}[grow=down]
  child {node [vertex,  label=180:{$1$}]{}  
  child {node [vertex, name=5,label=0:{$7$}]{}
  child {node [vertex, name=3,label=0:{$3$}]{} [sibling distance=1.8cm,level distance=15mm]
  child {node [vertex,name=2, label=-90:{$2$}]{}edge from parent[solid]node[yshift=3mm,xshift=-6mm]{$e_4$}}
  child {node [vertex, label=-90:{$5$}]{} node[yshift=2mm,xshift=-2mm]{$e_5$}}
  child {node [vertex, label=-90:{$4$}]{}  node[yshift=2mm,xshift=-4mm]{$e_6$}} }node[yshift=3mm,xshift=-3mm]{$e_7$}}}
  child {node [vertex, label=-90:{$6$}]{}edge from parent[solid] };
  \draw[red] (2)--(3);
  \draw[->, line width=1pt, shorten >=2pt, shorten <=2pt]  (2,-4) --(5.5,-4)
        node[yshift=4mm,xshift=-10mm] {$\phi_{e_4}$};
\end{tikzpicture}\quad
\begin{tikzpicture}
 [scale=0.5,vertex/.style={shape=circle, draw, inner sep=1.5pt, fill=black},
 subtree/.style={shape=ellipse, draw,minimum width=1.5cm, minimum height=.5cm},
 every fit/.style={ellipse,draw,inner sep=-2pt},
sibling distance=2cm,level distance=15mm,
 leaf/.style={label={[name=#1]below:$ $}},auto]
\node[vertex,name=1,label=0:{$0$}]{}[grow=down]
  child {node [vertex,  label=180:{$1$}]{}  
  child {node [vertex, label=0:{$7$}]{}
  child {node [vertex, label=0:{$2$}]{} [sibling distance=1.8cm,level distance=15mm]
  child {node [vertex, name=3, label=0:{$3$}]{}
  child {node [vertex,name=5, label=-90:{$5$}]{} node[yshift=3mm,xshift=-1mm]{$e_5$}}
  child {node [vertex, label=-90:{$4$}]{}  node[yshift=3mm,xshift=1mm]{$e_6$}} }}node[yshift=3mm,xshift=-3mm]{$e_7$}}}
  child {node [vertex, label=-90:{$6$}]{}edge from parent[solid] };
   \draw[red] (5)--(3);
  
\end{tikzpicture}\quad 

\begin{tikzpicture}
 [scale=0.6,vertex/.style={shape=circle, draw, inner sep=1.5pt, fill=black},
 subtree/.style={shape=ellipse, draw,minimum width=1.5cm, minimum height=.5cm},
 every fit/.style={ellipse,draw,inner sep=-2pt},
sibling distance=2cm,level distance=15mm,
 leaf/.style={label={[name=#1]below:$ $}},auto]
\node[vertex,name=1,label=0:{$0$}]{}[grow=down]
  child {node [vertex, label=180:{$1$}]{}  [sibling distance=12mm,level distance=12mm]
  child {node [vertex, name=5,label=0:{$7$}]{}
  child {node [vertex, label=0:{$2$}]{} 
  child {node [vertex, label=0:{$5$}]{}
  child {node [vertex, name=3, label=0:{$3$}]{} 
  child {node [vertex, name=4, label=0:{$4$}]{}  node[yshift=3mm,xshift=-3mm]{$e_6$} }} }}node[yshift=3mm,xshift=-3mm]{$e_7$}}}
  child {node [vertex, label=-90:{$6$}]{}edge from parent[solid]};
  \draw[red](3)--(4);
\draw[->, line width=1pt, shorten >=2pt, shorten <=2pt] (-5.5,-3) -- (-2.7,-3)
     node[yshift=4mm,xshift=-10mm] {$\phi_{e_5}$};
 
\end{tikzpicture} \quad
\begin{tikzpicture}
 [scale=0.6,vertex/.style={shape=circle, draw, inner sep=1.5pt, fill=black},
 subtree/.style={shape=ellipse, draw,minimum width=1.5cm, minimum height=.5cm},
 every fit/.style={ellipse,draw,inner sep=-2pt},
sibling distance=2cm,level distance=15mm,
 leaf/.style={label={[name=#1]below:$ $}},auto]
\node[vertex,name=1,label=0:{$0$}]{}[grow=down]
  child {node [vertex, name=3, label=180:{$1$}]{}  [sibling distance=12mm,level distance=12mm]
  child {node [vertex, name=5,label=0:{$7$}]{}
  child {node [vertex, label=0:{$2$}]{} 
  child {node [vertex, label=0:{$5$}]{}
  child {node [vertex, label=0:{$3$}]{} 
  child {node [vertex, label=0:{$4$}]{} }} }}node[yshift=3mm,xshift=-3mm]{$e_7$}}}
  child {node [vertex, label=-90:{$6$}]{}edge from parent[solid]};
  \draw[red](3)--(5);
\draw[->, line width=1pt, shorten >=5pt, shorten <=5pt] (-5.5,-3) -- (-2.5,-3)
     node[yshift=4mm,xshift=-10mm] {$\phi_{e_6}$};
 
\end{tikzpicture} 
\begin{tikzpicture}
 [scale=0.7,vertex/.style={shape=circle, draw, inner sep=1.5pt, fill=black},
 subtree/.style={shape=ellipse, draw,minimum width=1.5cm, minimum height=.5cm},
 every fit/.style={ellipse,draw,inner sep=-2pt},
sibling distance=2.2cm,level distance=15mm,
 leaf/.style={label={[name=#1]below:$ $}},auto]
\node[vertex,name=1,label=0:{$0$}]{}[grow=down]
  child {node [vertex, label=180:{$7(u_6)$}]{}  [sibling distance=14mm,level distance=14mm]
  child {node [vertex, name=5,label=-90:{$1(u_5)$}]{}node[xshift=-2mm,yshift=5mm]{\footnotesize$(x)$}}
  child {node [vertex, label=0:{$2(u_6)$}]{} [level distance=15mm]
  child {node [vertex, label=0:{$5$}]{}
  child {node [vertex, label=0:{$4$}]{} 
  child {node [vertex, label=-90:{$3(u_1)$}]{}  node[xshift=3mm,yshift=5mm]{\footnotesize$(x)$}}node[xshift=3mm,yshift=5mm]{\footnotesize$(x)$}} node[xshift=3mm,yshift=5mm]{\footnotesize$(y)$}}node[xshift=2mm,yshift=5mm]{\footnotesize$(x)$}}node[xshift=0mm,yshift=6mm]{\footnotesize$(y)$}}
  child {node [vertex, label=0:{$6(u_4)$}]{}edge from parent[solid] node[right]{\footnotesize$(y)$}};
  
\draw[->, line width=1pt, shorten >=5pt, shorten <=5pt] (-5,-4.5) -- (-2.5,-4.5)
     node[yshift=4mm,xshift=-8mm] {$\phi_{e_7}$};
 
\end{tikzpicture} 
\caption{An example of the bijection  $\Theta$. } \label{fig: prop1.3}
\end{figure}

We next prove Proposition~\ref{prop:G-sym-specialized}, which establishes further pairwise symmetries under the specialization \(x=y=1\).

\begin{proof}[Proof of Proposition~\ref{prop:G-sym-specialized}] 
Since the subgroup in Proposition~\ref{prop:G-sym-specialized} is generated by the  adjacent transpositions $(1\,2),(3\,4),(5\,6)$, it
suffices to prove the following three identities: 
    \begin{align}
G_n(1,1;u_1,u_2,u_3,u_4,u_5,u_6)&=G_n(1,1;u_2,u_1,u_3,u_4,u_5,u_6);\label{sym2} \\[5pt]
G_n(1,1;u_1,u_2,u_3,u_4,u_5,u_6)&=G_n(1,1;u_1,u_2,u_4,u_3,u_5,u_6);\label{sym3} \\[5pt]
G_n(1,1;u_1,u_2,u_3,u_4,u_5,u_6)&=G_n(1,1;u_1,u_2,u_3,u_4,u_6,u_5).\label{sym4} 
    \end{align}
It suffices to prove \eqref{sym2}; the proofs of the other two identities are analogous.  For an edge \(e\) of \(T\), let \(B_e\) and \(C_e\) be the two subtrees appearing in Fig.~\ref{treeABC}. Let \(\hat{E}(T)\) be the set of edges \(e\) for which \(B_e\) and \(C_e\) are of the types occurring in Case~(1) of Proposition~\ref{setsort}.  Define
\[
\hat{\Theta}(T)=\left(\prod_{e\in \hat{E}(T)}\phi_e\right)(T).
\]
Note that applying \(\phi_e\) leaves the corresponding edge in \(\hat{E}(T)\). Hence \(\hat{E}(\hat{\Theta}(T))=\hat{E}(T)\). Since the   maps \(\phi_e\) commute and each \(\phi_e\) is an involution, \(\hat{\Theta}\) is an involution on \(\mathcal P_n\).

For each edge \(e\in \hat{E}(T)\), Case (1) in Proposition~\ref{setsort}  shows that
\(\phi_e\) interchanges an elder twin leaf and a first singleton leaf.
All other relevant statistics are preserved. Hence
\[
{\rm fsleaf}(\hat{\Theta}(T))={\rm etleaf}(T),\qquad
{\rm etleaf}(\hat{\Theta}(T))={\rm fsleaf}(T),
\]
while
\[
{\rm ysleaf}(T),\quad {\rm entleaf}(T),\quad {\rm yerleaf}(T),\quad {\rm ryint}(T)
\]
are unchanged. Since \(x=y=1\), the possible changes in proper and improper
edges do not affect the weight. Summing over all \(T\in\mathcal P_n\), we
obtain \eqref{sym2}.

For \eqref{sym3} and \eqref{sym4}, one instead selects the edges corresponding to Case~(3) and Case~(2)  of Proposition~\ref{setsort}. In both cases, the selected set of edges is preserved by the corresponding partial product of   involutions, and the desired pair of statistics is interchanged while the remaining statistics are fixed.
\end{proof}

For example, Fig. \ref{fig: prop1.4} shows how \(\hat{\Theta}\) acts on a tree \(T\).   The weight changes from $u_2^2u_3u_6$ to $u_1^2u_3u_6$.
\begin{figure}[h]
 \begin{tikzpicture}
 [scale=0.7,vertex/.style={shape=circle, draw, inner sep=1.5pt, fill=black},
 subtree/.style={shape=ellipse, draw,minimum width=1.5cm, minimum height=.5cm},
 every fit/.style={ellipse,draw,inner sep=-2pt},
sibling distance=2cm, level distance=2cm,
 leaf/.style={label={[name=#1]below:$ $}},auto]
 
\node[vertex, label=90:{$0(u_6)$}]{}[grow=down]
  child {node [vertex, name=a, label=0:{$1$}]{} [sibling distance=10mm,level distance=15mm]
  child {node [vertex,name=b, label=-90:{$2$}]{} node[yshift=5mm,xshift=-1mm]{$e_1$} }
  child{node[vertex,label=-90:{$3(u_2)$}] {}  } }
  child {node [vertex, label=0:{$4$}]{}[sibling distance=12mm,level distance=15mm]
  child{node[vertex,label=-90:{$5(u_3)$}] {}   } }
  child {node [vertex, label=0:{$6$}]{}[sibling distance=10mm,level distance=15mm]
  child{node[vertex,label=-90:{$7$}] {} node[yshift=5mm,xshift=-1mm]{$e_2$}  } child{node[vertex,label=-90:{$8(u_2)$}] {}   } }  ;
   \draw[red, thick] (a) -- (b);
  \draw[->, line width=1pt, shorten >=2pt, shorten <=2pt]  (3,-1.7) --(5,-1.7)
    node[yshift=4mm,xshift=-7mm] {$\phi_{e_1}$};
\end{tikzpicture}
  \begin{tikzpicture}
 [scale=0.7,vertex/.style={shape=circle, draw, inner sep=1.5pt, fill=black},
 subtree/.style={shape=ellipse, draw,minimum width=1.5cm, minimum height=.5cm},
 every fit/.style={ellipse,draw,inner sep=-2pt},
sibling distance=1.7cm,level distance=2cm,
 leaf/.style={label={[name=#1]below:$ $}},auto]
 
\node[vertex, label=90:{$0$}]{}[grow=down]
  child {node [vertex, label=0:{$2$}]{} [sibling distance=12mm,level distance=15mm]
  child {node [vertex, label=0:{$1$}]{} 
  child{node[vertex,label=-90:{$3$}] {}}  } }
  child {node [vertex, label=0:{$4$}]{}[sibling distance=12mm,level distance=15mm]
  child{node[vertex,label=-90:{$5$}] {}   } }
  child {node [vertex,name=a, label=0:{$6$}]{}[sibling distance=10mm,level distance=15mm]
  child{node[vertex,name=b,label=-90:{$7$}] {} node[yshift=5mm,xshift=-1mm]{$e_2$}  } child{node[vertex,label=-90:{$8$}] {}   } };
   \draw[red, thick] (a) -- (b);
 
\end{tikzpicture}
   \begin{tikzpicture}
 [scale=0.7,vertex/.style={shape=circle, draw, inner sep=1.5pt, fill=black},
 subtree/.style={shape=ellipse, draw,minimum width=1.5cm, minimum height=.5cm},
 every fit/.style={ellipse,draw,inner sep=-2pt},
sibling distance=1.5cm,level distance=2cm,
 leaf/.style={label={[name=#1]below:$ $}},auto]
 
\node[vertex,name=5,label=90:{$0(u_6)$}]{}[grow=down]
  child {node [vertex, label=0:{$2$}]{} [sibling distance=12mm,level distance=15mm]
  child {node [vertex, label=0:{$1$}]{} 
  child{node[vertex,label=-90:{$3(u_1)$}] {}  } }}
  child {node [vertex, label=0:{$4$}]{}[sibling distance=12mm,level distance=15mm]
  child{node[vertex,label=-90:{$5(u_3)$}] {}   } }
  child {node [vertex, label=0:{$7$}]{}[sibling distance=12mm,level distance=15mm]
  child{node[vertex,label=0:{$6$}] {}   child{node[vertex,label=-90:{$8(u_1)$}] {}}   } }
  ;
   \draw[->, line width=1pt, shorten >=2pt, shorten <=2pt]  (-5,-3.3) --(-3,-3.3) node[yshift=4mm,xshift=-7mm] {$\phi_{e_2}$};
  
\end{tikzpicture}
  
\caption{An example of the bijection $\hat{\Theta}$.}
\label{fig: prop1.4}
\end{figure}

\noindent\textbf{Remark.} Case (4) in Proposition~\ref{setsort} shows that applying the   involution \(\phi_e\) to an edge whose associated subtrees \(B_e\) and \(C_e\) are of the same type preserves all the statistics. Therefore, the involutions used to prove \eqref{sym2}, \eqref{sym3}, and \eqref{sym4} are not unique: one may additionally compose them with any collection of such statistic-preserving involutions.

 For example, consider the rightmost plane tree shown in Fig. \ref{fig: prop1.4}. After applying \(\phi_{e_1}\) and \(\phi_{e_2}\), the resulting tree has the property that, for the edge \(e_3\), the associated subtrees \(B_e\) and \(C_e\) are of the same type, that is has more than one vertex. Hence, after   applying \(\phi_{e_3}\), it changes the shape of the tree, but preserves all statistics shown in Fig. \ref{fig: remark}. 
\begin{figure}[H]
\centering
 \begin{tikzpicture}
 [scale=0.7,vertex/.style={shape=circle, draw, inner sep=1.5pt, fill=black},
 subtree/.style={shape=ellipse, draw,minimum width=1.5cm, minimum height=.5cm},
 every fit/.style={ellipse,draw,inner sep=-2pt},
sibling distance=2cm, level distance=2cm,
 leaf/.style={label={[name=#1]below:$ $}},auto]
 
\node[vertex,name=5,label=90:{$0(u_6)$}]{}[grow=down]
  child {node [vertex, label=0:{$1$}]{} [sibling distance=10mm,level distance=15mm]
  child {node [vertex, label=-90:{$2$}]{} node[yshift=5mm,xshift=-1mm]{$e_1$} }
  child{node[vertex,label=-90:{$3(u_2)$}] {}  }node[yshift=10mm,xshift=7mm]{$e_3$} }
  child {node [vertex, label=0:{$4$}]{}[sibling distance=12mm,level distance=15mm]
  child{node[vertex,label=-90:{$5(u_3)$}] {}   } }
  child {node [vertex, label=0:{$6$}]{}[sibling distance=10mm,level distance=15mm]
  child{node[vertex,label=-90:{$7$}] {} node[yshift=5mm,xshift=-1mm]{$e_2$} } child{node[vertex,label=-90:{$8(u_2)$}] {}   } }
  ;
  \draw[->, line width=1pt, shorten >=2pt, shorten <=2pt]  (4,-1.7) --(8,-1.7)
    node[yshift=4mm,xshift=-15mm] {$\phi_{e_1}\phi_{e_2}\phi_{e_3}$};
\end{tikzpicture}
 \begin{tikzpicture}
 [scale=0.7,vertex/.style={shape=circle, draw, inner sep=1.5pt, fill=black},
 subtree/.style={shape=ellipse, draw,minimum width=1.5cm, minimum height=.5cm},
 every fit/.style={ellipse,draw,inner sep=-2pt},
sibling distance=2.5cm,level distance=1.5cm,
 leaf/.style={label={[name=#1]below:$ $}},auto]
 
\node[vertex,name=5,label=90:{$2$}]{}[grow=down]
  child {node [vertex, label=0:{$0(u_6)$}]{} [sibling distance=16mm,level distance=15mm]
  child {node [vertex, label=0:{$4$}]{}[sibling distance=12mm,level distance=15mm]
  child{node[vertex,label=-90:{$5(u_3)$}] {}   } }
  child {node [vertex, label=0:{$7$}]{}[sibling distance=12mm,level distance=15mm]
  child{node[vertex,label=0:{$6$}] {}   child{node[vertex,label=-90:{$8(u_1)$}] {}}   } }}
  child {node [vertex, label=0:{$1$}]{} 
  child{node[vertex,label=-90:{$3(u_1)$}] {}  } }
  ;
   
\end{tikzpicture}
\caption{An example illustrating the non-uniqueness of the involutions.}
\label{fig: remark}
\end{figure}

We now state three specializations of Propositions \ref{prop:G-sym} and \ref{prop:G-sym-specialized}. First, for $n \geq 1$, define
\begin{equation}
\begin{aligned}
\widetilde{G}_n & (x_{11}, x_{12}, x_2; y_{11}, y_{12}, y_2) \\
&= \sum_{T \in \mathcal{P}_n} x_{11}^{\text{sleaf}(T)} x_{12}^{\text{eleaf}(T)} x_2^{\text{yleaf}(T)} y_{11}^{\text{sint}(T)} y_{12}^{\text{eint}(T)} y_2^{\text{yint}(T)}
\end{aligned}  
\end{equation}
with the convention that $\widetilde{G}_0(x_{11}, x_{12}, x_2; y_{11}, y_{12}, y_2) = y_2$. For $T\in\mathcal{P}_n$, by definition, 
\begin{align}
{\rm sleaf}(T)&={\rm fsleaf}(T)+{\rm ysleaf}(T),\notag
\\[5pt]
{\rm eleaf}(T)&={\rm etleaf}(T)+{\rm entleaf}(T),\notag\\[5pt]
{\rm yleaf}(T)&={\rm secleaf}(T)+{\rm yerleaf}(T)
={\rm etleaf}(T)+{\rm yerleaf}(T),\label{ysecyer}
\end{align}
     and
\begin{align}
{\rm sint}(T)&={\rm fsleaf}(T)+{\rm ysleaf}(T),\notag\\[5pt]
{\rm eint}(T)&={\rm etleaf}(T)+{\rm entleaf}(T),\notag\\[5pt]
{\rm yint}(T)&={\rm gsint}(T)+{\rm ryint}(T) ={\rm fsleaf}(T)+{\rm ryint}(T).\label{ygry}
\end{align}
Therefore, 
\begin{align*}
\widetilde{G}_n(x_{11},x_{12},x_2;y_{11},y_{12},y_2)
=
G_n(1,1;
x_{11}y_{11}y_2,\,
x_{12}y_{12}x_2,\,
x_{11}y_{11},\,
x_{12}y_{12},\,
x_2,\,
y_2).
\end{align*}
Setting $x=y=1$ in  Propositions \ref{prop:G-sym}, we obtain 
\[
G_n(1,1;u_1,u_2,u_3,u_4,u_5,u_6)
=
G_n(1,1;u_2,u_1,u_4,u_3,u_6,u_5)
.\]
Then we obtain the following symmetry property.
\begin{cor}{\rm (\cite[Proposition 1.2]{Dong-Du-Ji-Zhang-2025}, \cite[Theorem 1.3]{Li-Lin-2025})} 
    For $n\geq 2$, we have 
\[
\widetilde{G}_n(x_{11},x_{12},x_2;y_{11},y_{12},y_2)
=
\widetilde{G}_n(x_{12},x_{11},y_2;y_{12},y_{11},x_2).
\]
\end{cor}

 \noindent\textbf{Remark.} Although the above corollary recovers the symmetry obtained by Li and Lin~\cite{Li-Lin-2025}, the involution used here is different from theirs. The following example compares the two maps on the same tree. The image under our involution $\Theta$ is shown in Fig.~\ref{Fig. Du-Ji-Zhang}, whereas the image under the involution of Li and Lin is shown in Fig.~\ref{Fig. Li-Lin}. The vertex types in Fig.~\ref{Fig. Du-Ji-Zhang} are read according to our right-to-left convention, while those in Fig.~\ref{Fig. Li-Lin} are read according to the left-to-right convention of Li and Lin.
 
 


\begin{figure}[H]
\centering
\begin{tikzpicture}[
    scale=0.7,
    every node/.style={font=\small},
    v/.style={circle, fill=black, inner sep=1.3pt},
    edge/.style={line width=0.45pt},
    arr/.style={->, line width=0.55pt}
]

\begin{scope}[shift={(0,0)}]
\node[v,label=above:{\(0(y_{12})\)}] (a0) at (0,0) {};
\node[v,label=left:{\(4(y_{12})\)}] (a1) at (-1,-1.1) {};
\node[v,label=below:{\(2(x_{12})\)}] (a2) at (1,-1.1) {};
\node[v,label=below:{\(3(x_{2})\)}] (a3) at (-2,-2.3) {};
\node[v,label=below:{\(1(x_{12})\)}] (a4) at (0,-2.3) {};

\draw(0.8,-0.5)node{$e_1$};
\draw(-0.8,-1.8)node{$e_2$};
\draw(-2,-1.8)node{$e_3$};
  \draw(-0.9,-0.5)node{$e_4$};
\draw[edge] (a0)--(a1);
\draw[red] (a0)--(a2);
\draw[edge] (a1)--(a3) (a1)--(a4);
\end{scope}

\node at (3,-1.2) {\(\phi_{e_1}\)};
\draw[arr] (2.2,-1.65)--(3.7,-1.65);

\begin{scope}[shift={(7,0)}]
\node[v,label=above:{\(2\)}] (a2) at (0,0) {};
\node[v,label=left:{\(4 \)}] (a1) at (-1,-1.1) {};
\node[v,label=below:{\(0 \)}] (a0) at (1,-1.1) {};
\node[v,label=below:{\(3 \)}] (a3) at (-2,-2.3) {};
\node[v,label=below:{\(1 \)}] (a4) at (0,-2.3) {};

\draw(-0.8,-1.8)node{$e_2$};
\draw(-2,-1.8)node{$e_3$};
  \draw(-0.9,-0.5)node{$e_4$};
\draw[edge] (a2)--(a1) (a0)--(a2);
\draw[edge] (a1)--(a3);
\draw[red](a1)--(a4);
\end{scope}

\node at (9.8,-1.2) {\(\phi_{e_2}\)};
\draw[arr] (9,-1.65)--(10.5,-1.65);

\begin{scope}[shift={(14,0)}]
\node[v,label=above:{\(2 \)}] (a2) at (0,0) {};
\node[v,label=left:{\(1 \)}] (a4) at (-1,-1.1) {};
\node[v,label=below:{\(0 \)}] (a0) at (1,-1.1) {};
\node[v,label=below:{\(3 \)}] (a3) at (-2,-2.3) {};
\node[v,label=below:{\(4 \)}] (a1) at (0,-2.3) {};

\draw(-2,-1.8)node{$e_3$};
  \draw(-0.9,-0.5)node{$e_4$};
\draw[edge] (a2)--(a4) (a0)--(a2);
\draw[edge]  (a1)--(a4);
\draw[red](a4)--(a3);
\end{scope}

\node at (2.5,-6) {\(\phi_{e_3}\)};
\draw[arr] (1.8,-6.5)--(3.5,-6.5);
\begin{scope}[shift={(6,-3)}]
\node[v,label=above:{\(2 \)}] (a2) at (0,-2) {};
\node[v,label=left:{\(3 \)}] (a3) at (-1,-3.1) {};
\node[v,label=below:{\(0 \)}] (a0) at (1,-3.1) {};
\node[v,label=left:{\(1 \)}] (a4) at (-1,-4.3) {};
\node[v,label=left:{\(4 \)}] (a1) at (-1,-5.3) {};
 \draw(-0.8,-2.5)node{$e_4$};
  \draw[red](a2)--(a3);
\draw[edge]  (a0)--(a2);
\draw[edge] (a4)--(a3) (a1)--(a4);
\node at (3,-3.1) {\(\phi_{e_4}\)};
 \draw[arr] (2.2,-3.5)--(3.8,-3.5);
\end{scope}
\begin{scope}[shift={(13.5,-3)}]
\node[v,label=above:{\(3(y_{2})\)}] (a3) at (0,-2) {};
\node[v,label=left:{\(2(y_{11})\)}] (a2) at (-1,-3.1) {};
\node[v,label=below:{\(0(x_{11})\)}] (a0) at (-1,-4.1) {};
\node[v,label=right:{\(1(y_{11})\)}] (a4) at (1,-3.1) {};
\node[v,label=below:{\(4(x_{11})\)}] (a1) at (1,-4.1) {};

\draw[edge] (a2)--(a3) (a0)--(a2);
\draw[edge] (a4)--(a3) (a1)--(a4);

\end{scope}

\end{tikzpicture}
\caption{An example of the involution
\(\Theta\).}
\label{Fig. Du-Ji-Zhang}
\end{figure}

\begin{figure}[H]
\centering
\begin{tikzpicture}[
    scale=0.7,
    every node/.style={font=\small},
    v/.style={circle, fill=black, inner sep=1.3pt},
    edge/.style={line width=0.45pt},
    arr/.style={->, line width=0.55pt}
]

\begin{scope}[shift={(0,0)}]
\node[v,label=above:{\(0(y_{2})\)}] (a0) at (0,0) {};
\node[v,label=left:{\(4(y_{12})\)}] (a1) at (-1,-1.1) {};
\node[v,label=below:{\(2(x_{12})\)}] (a2) at (1,-1.1) {};
\node[v,label=below:{\(3(x_{12})\)}] (a3) at (-2,-2.3) {};
\node[v,label=below:{\(1(x_2)\)}] (a4) at (0,-2.3) {};

\draw[edge] (a0)--(a1) (a0)--(a2);
\draw[edge] (a1)--(a3) (a1)--(a4);
\end{scope}

\node at (3,-1.2) {\(\rho\)};
\draw[arr] (2.2,-1.65)--(3.7,-1.65);

\begin{scope}[shift={(5.5,0)}]
\node[v,label=above:{\(2\)}] (b2) at (0,0) {};
\node[v,label=left:{\(4\)}] (b1) at (-0.8,-1.0) {};
\node[v,label=right:{\(1\)}] (b4) at (0,-2.0) {};
\node[v,label=left:{\(3\)}] (b3) at (-0.8,-3.0) {};

\draw[edge]  (b2)--(b1);
\draw[edge] (b1)--(b4) (b4)--(b3);
\end{scope}

\node at (7.8,-1.2) {\(\psi\)};
\draw[arr] (7,-1.65)--(8.5,-1.65);

\begin{scope}[shift={(9.3,0)}]
\node[v,label=above:{\(2\)}] (c2) at (0,0) {};
\node[v,label=right:{\(4\)}] (c1) at (0.8,-1.0) {};
\node[v,label=left:{\(1\)}] (c4) at (0,-2.0) {};
\node[v,label=right:{\(3\)}] (c3) at (0.8,-3.0) {};

\draw[edge] (c1)--(c4) (c4)--(c3);
\draw[edge] (c1)--(c2);
\end{scope}

\node at (12.3,-1.2) {\(\rho^{-1}\)};
\draw[arr] (11.4,-1.65)--(13,-1.65);

\begin{scope}[shift={(15.5,0)}]
\node[v,label=above:{\(0(y_2)\)}] (d0) at (0,0) {};
\node[v,label=right:{\(2(y_2)\)}] (d2) at (0,-1) {};
\node[v,label=right:{\(4(x_2)\)}] (d1) at (1.0,-2) {};
\node[v,label=right:{\(1(y_{11})\)}] (d4) at (-1.0,-2) {};
 
\node[v,label=below:{\(3(x_{11})\)}] (d3) at (-1.0,-3) {};

\draw[edge] (d0)--(d2) (d2)--(d1);
\draw[edge] (d4)--(d2) (d3)--(d4);
\end{scope}

\end{tikzpicture}
\caption{An example of the involution
\(\widetilde{\Theta}=\rho^{-1}\circ\psi\circ\rho\).}
\label{Fig. Li-Lin}
\end{figure}

Let \(\mathcal P_n^{\rm ta}\) be the subset of \(\mathcal P_n\) consisting of trees with no young interior vertices. Equivalently, under our right-to-left convention, every interior vertex has a rightmost child that is a leaf. Up to reflection, these are the tip-augmented plane trees introduced by Donaghey~\cite{Donaghey-1977}. 

A plane tree belongs to  \(\mathcal P_n^{\rm ta}\) if and only if
\[
{\rm gsint}(T)={\rm ryint}(T)=0.
\]
Since \({\rm gsint}(T)={\rm fsleaf}(T)\), this is equivalent to
\[
{\rm fsleaf}(T)={\rm ryint}(T)=0.
\]
For $n \geq 1$,  define
\begin{equation}
\begin{aligned}
M_n & (v_{1}, v_{2}, v_3; w_{1}, w_{2}) \\
&= \sum_{T \in \mathcal P^{\rm ta}_n} v_{1}^{\text{sleaf}(T)} v_{2}^{\text{etleaf}(T)} v_3^{\text{entleaf}(T)} w_{1}^{\text{yerleaf}(T)}  w_2^{\text{secleaf}(T)}
\end{aligned}  
\end{equation}
with the convention that $M_0(v_{1},  v_2,v_3; w_{1}, w_{2}) = v_2$.

Since each elder twin leaf is paired with a unique second leaf, we have 
\begin{equation} \label{rel:Mot}
M_n(v_1,v_{2}, v_3; w_{1}, w_{2})=G_n(1,1;0,v_2w_2,v_1,v_3,w_1,0).
\end{equation}

 Applying \eqref{sym3} into \eqref{rel:Mot}, we obtain the following corollary.
\begin{cor}{ \rm (\cite[Proposition 1.3]{Dong-Du-Ji-Zhang-2025},   \cite[Proposition 1.8]{Li-Lin-2025})}
For \(n\geq 2\), we have
\[
M_n(v_1,v_{2}, v_3; w_{1}, w_{2})
=
M_n(v_3,v_{2}, v_1; w_{1}, w_{2}).
\]
\end{cor}

Let \(\mathcal P_n^{\rm ns}\) be the subset of  \(\mathcal P_n\) consisting of trees with  no singleton leaves.  The condition
\(T\in\mathcal P_n^{\rm ns}\) is equivalent to
\[
{\rm fsleaf}(T)={\rm ysleaf}(T)=0.
\]
Thus this subclass is extracted from \(G_n\) by setting \(u_1=u_3=0\). Applying \eqref{sym4}, we get that for \(n\geq 2\), 
\[
G_n(1,1;0,u_2,0,u_4,u_5,u_6)
=
G_n(1,1;0,u_2,0,u_4,u_6,u_5).
\]

 
  Under the standard depth-first,  encoding of plane trees by Dyck paths, leaves
correspond to peaks \(UD\). For \(n\ge2\), singleton leaves correspond exactly
to consecutive occurrences of \(UUDD\). Thus \(\mathcal P_n^{\rm ns}\)
corresponds to Dyck paths avoiding \(UUDD\), whose enumeration is recorded in OEIS~A082582 \cite{OEIS}.


 \vskip 0.2cm
\noindent{\bf Acknowledgment.} 
The authors thank the referees for their insightful comments and suggestions. This work was supported by the National Natural Science Foundation of China.


\begin{thebibliography}{0}
	\setlength{\itemsep}{-.8mm}
	\addcontentsline{toc}{section}{}

\bibitem{Archer-Gregory-Pennington-Slayden-2019} K. Archer, A. Gregory, B. Pennington, and S. Slayden, Pattern restricted quasi-Stirling permutations, Australas. J. Combin., 74 (2019), 389. 

\bibitem{Chapoton-2002} F. Chapoton, Op\'erades diff\'erentielles gradu\'ees sur les simplexes et les permuto\`edres, Bull. Soc. Math. France, 130 (2002), 233--251.

\bibitem{Chen-Deutsch-Elizalde-2006} W. Y. C. Chen, E. Deutsch, and S. Elizalde, Old and young leaves on plane trees, European J. Combin., 27 (2006), 414--427.

\bibitem{Chen-Fu-Wang-2025}
W.~Y.~C.~Chen, A.~M.~Fu and E.~L.~Wang, A grammatical calculus for the Ramanujan polynomials,  arXiv:2506.01649v1.

\bibitem{Chen-Yang-2021}
W.~Y.~C.~Chen and H.~R.~L.~Yang,  A context-free grammar for the Ramanujan-Shor polynomials,
Adv. in Appl. Math., 126 (2021), Paper No. 101908, 24.

\bibitem{Donaghey-1977} R. Donaghey, Restricted plane tree representations of four Motzkin-Catalan equations, J. Combin. Theory Ser. B, 22 (1977), 114--121.

\bibitem{Dong-Du-Ji-Zhang-2025} J.~J.~W. Dong, L.~R. Du, K.~Q. Ji and  D.~T.~X. Zhang, New refinements of Narayana polynomials and Motzkin polynomials, Adv. in Appl. Math., 166  (2025), Paper No. 102855, 40. 


  
\bibitem{Dumont-Ramamonjisoa-1996}
D. Dumont and A. Ramamonjisoa, Grammaire de Ramanujan et arbres de Cayley,
Vol.~3, 1996, Research Paper 17, approx.~18, the Foata Festschrift.

  \bibitem{Elizalde-2021} S. Elizalde,  Descents on quasi-Stirling permutations,  J. Combin. Theory Ser. A, 180 (2021),  Paper No. 105429, 35.

\bibitem{Gessel-Stanley-1978}
I.~Gessel and R.~P.~Stanley, Stirling polynomials, J. Combin. Theory Ser. A, 24 (1978), 24--33.

\bibitem{Guo-Zeng-2007}
V.~J.~W.~Guo and J.~Zeng, A generalization of the Ramanujan polynomials and plane trees,
 Adv. in Appl. Math., 39 (2007), 96--115.

\bibitem{Janson-2008}
S.~Janson, Plane recursive trees, Stirling permutations and an urn model,
 in:  Fifth Colloquium on Mathematics and Computer Science, vol.~AI of Discrete Math. Theor. Comput. Sci. Proc., Assoc. Discrete Math. Theor. Comput. Sci., Nancy 2008, 541--547.

\bibitem{Janson-Kuba-Panholzer-2011}
S.~Janson, M.~Kuba and A.~Panholzer, Generalized Stirling permutations, families of increasing trees and urn models, J. Combin. Theory Ser. A, 118 (2011), 94--114.

\bibitem{Koganov-1996} L. M. Koganov, Universal bijection between Gessel-Stanley permutations and diagrams of connections of corresponding ranks, Uspekhi. Mat. Nauk, 51  (1996), 165-166.

\bibitem{Kuba-Panholzer-2007}
  M. Kuba and A. Panholzer, On the degree distribution of the nodes in increasing trees, J. Combin. Theory Ser. A, 114 (2007), 597--618.

\bibitem{Li-Lin-2025} Y. Li and Z. Lin, Bijections in weakly increasing trees via binary trees, Discrete Math., 348 (2025), Paper No. 114632, 17. 

\bibitem{OEIS}
OEIS Foundation Inc., The On-Line Encyclopedia of Integer Sequences.

\bibitem{Remmel-Wilson-2015}
J.~B.~Remmel and A.~T.~Wilson, Block patterns in Stirling permutations, J. Comb., 6 (2015), 179--204.

\bibitem{Shor-1995}
P.~W.~Shor, A new proof of Cayley's formula for counting labeled trees, J. Combin. Theory Ser. A, 71 (1995), 154--158.

\bibitem{Stanley-2015} R.~P.~Stanley, Catalan Numbers, Cambridge University Press, 2015. 

\bibitem{Stanley-2024}
R.~P.~Stanley, Enumerative Combinatorics. Vol.~2, vol.~208 of Cambridge Studies in Advanced Mathematics, Cambridge University Press, Cambridge[2024], 2nd ed., with an appendix by Sergey Fomin.

\bibitem{Zeng-1999}
J.~Zeng, A Ramanujan sequence that refines the Cayley formula for trees, Ramanujan J., 3 (1999), 45--54.

\end{thebibliography}

\end{document}